\documentstyle[12pt]{article}
\begin{document}
\addtolength{\baselineskip}{6pt}
\title{The theta divisor of the bidegree (2,2) threefold
in ${\bf P}^{2} \times {\bf P}^{2}$ }
\author{Atanas Iliev\footnotemark[1]}
\date{ }
\maketitle
\footnotetext{Supported in part by the Bulgarian foundation
"Scientific Research" and
by NSF under the US-Bulgarian project
"Algebra and Algebraic Geometry".}

\centerline{\sc 0. Introduction}
\smallskip

{\bf (0.1).}
In this paper we apply a new approach to the study of the theta
divisor of a standard conic bundle. As an example we examine the
Verra threefold $T = T(2,2)$ -- the divisor of bidegree (2,2) in
${\bf P}^{2}\times{\bf P}^{2}$.
 The threefold $T$ deserves a special
attention because of the recent observation of A.Verra (see [Ve])
that the existence of two conic bundle structures on $T$ implies
a new counterexample to the Torelli theorem for Prym varieties.
Moreover, this counterexample is not related to the 4-gonal
correspondence (see [Do]) of Donagi
(which had covered all the known
non-trivial counterexamples).

\smallskip
{\bf (0.2).}
Let  $p: T \longrightarrow {\bf P}^{2}$ be any of the two
 projections
which make  $T$ a standard conic bundle over ${\bf P}^{2}$.
 In section 3 we
show that the existence of $p$ implies the existence of a
 special family
${\cal C}_{\theta}$ of curves on $T$, which is mapped
 (via the Abel-Jacobi
map) onto a copy of the theta divisor ${\Theta} (T)$
 (see Theorem 4.1). The
general curve $C \in {\cal C}_{\theta}$ is an elliptic curve of
bidegree (3,6) which lies in a hyperplane section of $T$. By
construction, the curves of ${\cal C}_{\theta}$ parameterize the
minimal sections (resp. -- the maximal subbundles of rank 1)
 of a special
family of rank 2 vector bundles (resp. -- ruled surfaces) over the
space of plane cubics (see section 2, (3.1), [LN], [Se]).
 It turns out
that ${\cal C}_{\theta}$ is generically a 2-sheeted covering of the
family of effective divisors $Supp \ {\Theta}$ related to the
 Wirtinger
description of the theta divisor (see sect. 2,3,4, and (1.2.2)).
 The last is
used in the proof of Theorem 4.1: The Abel-Jacobi image of
 ${\cal C}_{\theta}$
is a copy of $\Theta(T)$. As a second application of the connection
 between
${\cal C}_{\theta}$ and  $\Theta$ we prove that the general
 hyperplane section
of $T(2,2)$ is a K3 surface, which is a double covering of
 ${\bf P}^{2}$
branched over a sextic with a vanishing theta null
 (see (5.2.1) and (5.3)).

\smallskip
{\bf (0.3).}
In section 6, we separate a 3-dimensional component $Z  \subset
Sing \ {\Theta}$ as an Abel-Jacobi image of the 6-dimensional
 family of
elliptic
sextics of bidegree (3,3) on $T$. The connection between the
geometry
of $T$ and $Sing \ {\Theta}$ , based on the study of the last
 family of
curves,
is used to  prove the Torelli theorem for $T(2,2)$:
 The threefold $T \subset
{\bf P}^{8}$ coincides with the intersection of the tangent
 cones of $\Theta$
in the points of $Z$ (see Theorem 6.6). Note that the Torelli
 theorem, stated
in this form, is not a direct consequence of general facts
 about Prym
varieties: one has also to see that the projectivized tangent cone
$Cone_{z}$ ($z$ -- a general point of $Z$) does not belong to the
 ``trivial''
component $D_{6}(W)$ of the determinantal locus $D_{6}(T)$ -- see
(6.3.3), (6.5.1), and the proof of Theorem (6.6).

In his paper [Ve], A.Verra proves that the two discriminant
pairs, which correspond to the two conic bundle structures on
 $T(2,2)$, fits
into the classical Dixon construction. In section 7 we show that,
 in case
of a nodal $T(2,2)$, this Dixon correspondence can be represented
 as a
composition of two 4-gonal correspondences of Donagi
 (see Corollary 7.3).

\bigskip
\centerline{\sc 1. Preliminaries}
\smallskip

{\bf (1.1).}
{\it The bidegree (2,2) divisor T} (see [Ve]).
\smallskip

{\bf (1.1.1).}
Let $seg: {\bf P}^{2}\times{\bf P}^{2} \longrightarrow{\bf P}^{8}$
 be the
Segre embedding, and let $W \subset {\bf P}^{8}$ be the image of
 $seg$. Let
$p_{1}:W \rightarrow {\bf P}^{2}$ and
 $p_{2}:W \rightarrow {\bf P}^{2}$
be the canonical projections. Denote
 by ${\cal O}_{W}(m,n)$ the sheaf
$p_{1}^{*}{\cal O}_{{\bf P}^{2}}(m)\otimes p_{2}^{*}{\cal O}_{{\bf
P}^{2}}(n)$,
m,n -- integers.

{\bf (1.1.2).}
Let $T = T(2,2) \in \mid{\cal O}_{W}(2,2)\mid $, be a
{\it bidegree (2,2) divisor} on $W$, and let
$p_{i}:T \rightarrow {\bf P}^{2}$ be the restrictions of
 $p_{i}$ on $T$.

{\bf (1.1.3).}
Let, moreover, $T = T(2,2)$ be {\it smooth}.  Then $p_{i}$
 defines a
standard conic bundle structure on $T$, $i = 1,2$.  Let $i \in
\{1,2\}$ be fixed, and let ${\Delta}_{i} = \{ x \in {\bf P}^{2} :
sing \ p^{-1}(x) \neq \oslash \} = \{ x \in {\bf P}^{2} :
 p^{-1}(x) = l
+ \overline{l}$ is a plane conic of rank 2 \} be the discriminant
curve of $p_{i}$.  It is not hard to see that ${\Delta}_{i}$ is a
smooth plane sextic ( $T$ is smooth).  Denote by

$\tilde{\Delta}_{i}$ = $ \{ l - \ \mbox{a line in}\  T :
  \exists x \in
{\Delta},\ s.t.\  l \ \mbox{is a component of}\ {p}^{-1}(x) \} $
 the double
discriminant curve of $p_{i}$. With a probable abuse of the
 notation, we
denote by
$p_{i}: \tilde{\Delta}_{i} \rightarrow {\Delta}_{i}$
also the induced (unbranched) double covering. Let $\eta_{i} \in
Pic^{o}{\Delta}_{i} $ be the torsion-2 sheaf, which defines
 $p_{i}.$

{\bf (1.1.4).}
Let $J(T)$ be the intermediate jacobian of $T$, and let
$P_{i} = P({\Delta}_{i},{\eta}_{i})$ be the Prym variety of
$({\Delta}_{i},{\eta}_{i})$. It is well-known, that $J(T)$
 and $P_{i}$
are  isomorphic  as   principally   polarized   abelian
 varieties
(see[B1]).

In particular, $dim J(T) = dim P_{i} = 9$. It follows immediately
 that
$P({\Delta}_{1},{\eta}_{1})$ and  $P({\Delta}_{2},{\eta}_{2})$
are isomorphic as p.p.a.v..

{\bf (1.1.5).}
In [Ve], A.Verra proves that the discriminant pairs $({\Delta}_{1},
{\eta}_{1})$ and $({\Delta}_{2},{\eta}_{2})$ correspond to each
 other
by the classical Dixon construction.
  Moreover, let ${\cal P}_{6} = \{
(\Delta , \eta ) :  \Delta$ is a smooth plane curve of degree 6,
 and
$\eta$ is a torsion-2 sheaf in ${\bf Pic}({\Delta}) \}$, and let
$p_{6} :  {\cal P}_{6} \longrightarrow {\cal A}_{9}$ be the Prym
map.
Then  $deg \ p_{6} = 2$, and $p_{6}$ is branched along the locus
of intermediate jacobians of nodal quartic double solids,
 see [Ve]. The
general fiber of $p_{6}$ consists of a couple of pairs
$(({\Delta}_{1},{\eta}_{1}),({\Delta}_{2},{\eta}_{2}))$,
 which arises
from a bidegree (2,2) divisor $T$.

{\bf (1.1.6).}
We call the smooth bidegree (2,2) divisor $T \in W$
{\it the Verra threefold}.

\smallskip
{\bf (1.2).}
{\it The intermediate jacobian $J(T)$ as a Prym variety.}
\smallskip

{\bf (1.2.1.)}
Let $\Theta$ be the divisor of the principal polarization
 (the theta
divisor) of $J(T)$. It follows from the preceding that we can
identify
$\Theta$ and the theta divisor of the Prym variety $P_{i}$,
 i = 1,2.
Then, because of the Wirtinger description of Prym varieties,
we can
describe $J = J(T), \Theta$, etc., only in terms
of $\tilde{\Delta}_{i}$ and ${\Delta}_{i}$. As a direct
 corollary we
obtain:

{\bf (1.2.2).}
Let $p_{i}$ be fixed. Then:

The jacobian $J(T)$ is isomorphic to
 $P(\tilde{\Delta}_{i},{\Delta}_{i}) =$

$ \{ {\cal L} \in {\bf Pic}^{18}\tilde{\Delta}_{i}: Nm{\cal L} =
{\omega}_{{\Delta}_{i}},
\ \mbox{and} \   h^{0}({\cal L}) \ \mbox{even} \ \} $.

${\Theta}(T) \cong
{\Theta}_{i} = \{ {\cal L} \in J(T) : h^{0}({\cal L}) \ge 2 \}$.

There exists a subset in $Sing \ {\Theta}$, which is isomorphic
 to the
set of stable singularities of ${\Theta}$, with respect to
 $p_{i}$, i.e.:
$$
Sing^{st}_{i}(\Theta) =
 \{ {\cal L} \in {\Theta}: h^{0}({\cal L}) \ge 4 \};
$$
similarly -- for the exceptional singularities of
 $\Theta$, w.r. to $p_{i}$.

In section 6, we shall describe a component
 $Z$ of $Sing \ {\Theta}$, the
points of which are stable w.r.to both  $p_{1}$ and $p_{2}$.

\smallskip
{\bf (1.3).}
{\it Minimal sections of ruled surfaces over elliptic curves}
\smallskip

{\bf (1.3.0).}
Here we collect some facts about ruled surfaces
 (esp. -- on elliptic curves),
which will be used in section 2, see [H, ch.V.2], [LN], [Se].

{\bf (1.3.1).}
Let $C$ be a smooth curve.  By definition, a ruled surface $p:S
\rightarrow C$ is a surface which can be represented
 in the form $ S =
{\bf P}_{C}({\cal E}) $, where ${\cal E}$ is a locally free
sheaf of
rank 2 (a rank 2 vector bundle) over $C$.  The representation $S =
{\bf P}_{C}({\cal E})$ is unique up to a multiplication by an
invertible sheaf:
  ${\bf P}_{C}({\cal E}) \cong {\bf P}_{C}({\cal E}
\otimes {\cal L})$, ${\cal L} \in {\bf Pic}(C)$.

The bundle ${\cal E}$ is called {\it normalized} if
$h^{0}({\cal E}) \neq 0$,  but
$h^{0}({\cal E} \otimes {\cal L}) = 0$
for any invertible ${\cal L}$ such that $deg{\cal L} < 0$.

{\bf (1.3.2).}
Obviously, any ruled surface $S$ has a representation
$S = {\bf P}_{C}({\cal E})$,
for some normalized ${\cal E}$. Anyway, such a representation is,
in general, far from unique (see e.g. [LN, Cor.3.2],
 or (1.3.4) -- in case
$g(C) = 1$).

Let ${\cal E}$ be normalized, and let
$C_{o} \in \mid {\cal O}_{{\bf P}({\cal E})/C}(1) \mid $
be a tautological section for ${\cal E}$. The invariant property of
such a $C_{o}$ is that $C_{o}$ is a section of $p:S \rightarrow C$,
for which the number
$-e(C) = (C.C)_{S},  C - \ \mbox{a section of} \  S$
is minimal, i.e. $C_{o}$ is a {\it minimal section} of $S$.
 The  number
$e = e(S) = -(C_{o}.C_{o})$
is an {\it invariant} of $S$.

{\bf (1.3.3).}
The surface
$S = {\bf P}({\cal E})$
is called {\it decomposable} if ${\cal E}$ is decomposable (see e.g.
[H,ch.V.2]). Otherwise, $S$ is called {\it indecomposable} (ibid).

{\bf (1.3.4.).}
The cardinality of the set of minimal sections of $S$
 closely depends
on the decomposability of $S$, and on (the parity) of
 the invariant
$e = e(S)$  (see [LN, Cor.3.2]).

In section 3, we shall use the description of these sets
 only in case
$g(C) = 1$. All of the following can be found in [H, ch.V.2]:

\smallskip
{\bf (*)}.
{\it Minimal sections of a ruled surface over an elliptic curve.}
\smallskip

Let $C$ be an elliptic curve, and let $S$ be a ruled surface
 on $C$.
Then, one of the following alternatives is valid:

(1). $S$ is decomposable, $e = e(S) > 0$. The normalized sheaf
 ${\cal E}$
     for $S$ is unique, and the minimal section $C_{o}$ is unique.

(2). $S$ is decomposable, $e = e(S) = 0$. In this case
     $S = {\bf P}_{C}({\cal O} \oplus {\varepsilon})$,
     $deg({\varepsilon}) = 0$, and {\it either}:

(a). ${\varepsilon} = {\cal O}_{C}, \ S = C \times {\bf P}^{1}$,
     and the set of minimal sections of $S$ is parameterized by
     the points of ${\bf P}^{1}$, {\it or}:

(b). ${\varepsilon} \neq {\cal O}_{C}$. Then the normalized sheaf
     ${\cal E}$ can be chosen in two ways:
     ${\cal E}^{+} = {\cal E}$ and
     ${\cal E}^{-} = {\bf P}( {\cal O} \oplus (-{\varepsilon}) )$,
     where
     $-{\varepsilon} = {\varepsilon}^{\otimes -1}$.
     Correspondingly, there are exactly two minimal sections of $S$:
     $C^{+}$ and $C^{-}$ (each -- the unique tautological section
     of the corresponding normalized bundle).

(3). $S$ is the unique indecomposable ruled surface, s.t. $e(S) = 0$.
     Then the normalized sheaf for $S$ is unique, and the
     corresponding
     minimal section is unique.

(4). $S$ is the unique indecomposable ruled surface, s.t.
     $e(S) = -1$.
     Then, the set

     $\{ {\cal E} - normalized : S = {\bf P}({\cal E}) \}$
     is parameterized by the points of the elliptic curve $C$.
     For any such
     ${\cal E}$, the minimal section $C_{o}({\cal E})$ is unique.

\bigskip
\centerline{\sc 2. Minimal sections of the canonical conic
 bundle surfaces}
\smallskip

{\bf (2.0).}
Everywhere in this section the conic bundle structure
$p:T \rightarrow {\bf P}^{2}$
is fixed; we let
$p = p_{1}, {\Delta} = {\Delta}_{1}, {\eta} = {\eta}_{1}$, etc.

\smallskip
{\bf (2.1).}
{\it The sets $ Supp \ {\Theta}$ and $Supp \ P^{-}$.}
\smallskip

{\bf (2.1.1).}
Let $p = p_{1}: T \rightarrow {\bf P}^{2}$, etc., be as above.
Let $Nm:{\bf Pic}^{18}{\tilde{\Delta}} \rightarrow
 {\bf Pic}^{18}{\Delta}$
be the norm map (see [ACGH, app.C]) Then
 $Nm^{-1}({\omega}_{\Delta})$ splits
into two components:

$P^{+} = \{ {\cal L} \in Nm^{-1}({\omega}_{\Delta}) : h^{0}({\cal L})
\ \mbox{even}  \}$ and
$P^{-} =$ the same, but $h^{0}({\cal L}) \ odd.$

\smallskip
{\bf (2.1.2).}
{\it $Supp\ {\Theta}$ and $Supp(P^{-}).$ }

The general sheaf ${\cal L} \in P^{+}$ is non-effective, i.e. the
linear
system of effective divisors $\mid {\cal L} \mid$ is empty. However,
 the
subset of the effective sheaves ${\cal L} \in P^{+}$ is exactly
 the theta
divisor ${\Theta}$ (see (1.2.2)). This gives a reason to define
 the set

$Supp \ {\Theta} :=
 \{ L \in \mid {\cal L} \mid : {\cal L} \in {\Theta} \}$,

i.e., $Supp \ {\Theta}$ is the set of all effective
 divisors in the linear
systems of the sheaves ${\cal L} \in {\Theta}$. Similarly:

$Supp(P^{-}) := \{ L \in \mid {\cal L} \mid : {\cal L} \in P^{-}\}$
(all the sheaves ${\cal L} \in P^{-}$ are effective).

\smallskip
{\bf (2.1.3).}
{\it The maps $p_{*}^{+}: Supp \ {\Theta} \rightarrow
\mid {\cal O}_{{\bf P}^{2}}(3) \mid$ and
$p_{*}^{-}: Supp(P^{-}) \rightarrow
 \mid {\cal O}_{{\bf P}^{2}}(3) \mid$.}

The set ${\Theta} \cup P^{-}$ coincides with the set of all
 ``effective''
sheaves in the preimage
$Nm^{-1}({\omega}_{\Delta})$. Moreover, on the level of
 effective divisors,
the map $Nm$ coincides with the usual projection
$p_{*}: Symm^{18}{\tilde {\Delta}} \rightarrow Symm^{18}{\Delta}.$
In particular, if ${\cal L} \in Nm^{-1}{\omega}$ is effective and
$L \in \mid {\cal L} \mid$, then
${\cal O}(p_{*}L) = Nm{\cal L} = {\omega}_{\Delta}
 = {\cal O}_{\Delta}(3).$
Since $deg\ {\Delta} = 6$, the linear system
 $\mid {\cal O}_{\Delta}(3) \mid$
is isomorphic to $\mid {\cal O}_{{\bf P}^{2}}(3) \mid$.
 In particular, the
effective divisor $p_{*}L \in Symm^{18}{\Delta}$ is a scheme
 intersection
of ${\Delta}$ and (a unique) plane cubic curve $C(L) = p_{*}(L).$
In particular, after  composing with the corresponding
 restriction maps, the
map $p_{*}$ defines the maps:

$ p_{*}^{+}: Supp \ {\Theta}  \rightarrow
 \mid {\cal O}_{\bf P^2}(3)\mid $
(= the set of all the plane cubics) and

$p_{*}^{-}: Supp \ (P^{-}) \rightarrow
\mid{\cal O}_{{\bf P}^{2}}(3)\mid .$

It is not hard to see (see e.g. [Sh, Lemma 3.20]) that:

(1). the maps $p_{*}^{+}$ and $p_{*}^{-}$ are surjective;

(2). the general fibers of $p_{*}^{+}$ and $p_{*}^{-}$ are finite.

In (2.2), we shall describe the general fibers of these two maps.

\smallskip
{\bf (2.2).}
{\it The canonical conic bundle surface $S(C)$ and the preimage
$p_{*}^{-1}(C)$.}
\smallskip

{\bf (2.2.1).}
Let $C$ be a sufficiently general plane cubic curve. In particular,
$C$ can be supposed to be smooth, and intersecting the discriminant
sextic ${\Delta}$ in 18 disjoint points $x_{1}, ..., x_{18}.$

The non-minimal ruled surface $S(C) = p^{-1}(C) \subset T = T(2,2)$
is a standard conic bundle over the (elliptic) curve $C$;
 the degenerate
fibers of $p:S(C) \rightarrow C$ are
  $f_{i} = p^{-1}(x_{i}) = l_{i} +
{\overline{l}_{i}}, i=1,...,18.$

We call a surface $p^{-1}(C) \subset T$ , $C$ -- any plane cubic
(resp.$C$ -- a general plane cubic), a {\it canonical  conic  bundle
surface}
on $T$ (resp. -- a {\it general c.c.b.s.}) -- w.r.to $p = p_{1}.$

\smallskip
{\bf (2.2.2).}
{\it The set ${\Sigma}(C)$.}

Let the cubic $C$ be as in (2.2.1), and let

${\Sigma}(C) =
\{ {\sigma: \bigcup \{ x_{i} \} \rightarrow
\bigcup \{ l_{i}, \overline {l}_{i} \} :
 \{\sigma}(x_{i}) \in \{ l_{i}, \overline{l}_{i} \}, i=1,...,18 \}$

be the set of ``choice'' maps for $C$.
We can define, in an obvious way, the two-argument signature
 function:

$sgn: {\Sigma}(C) \times {\Sigma}(C) \longrightarrow \{ +1,-1 \}$
as follows:

$sgn ({\sigma}',{\sigma}'') = +1$, if
$\# (Image({\sigma}')\cap Image({\sigma}'')) \in  2{\bf Z}$,

otherwise
$sgn({\sigma}',{\sigma}'') = -1.$

\smallskip
{\bf (2.2.3).}
{\it The map $L:{\Sigma}(C) \rightarrow p_{*}^{-1}(C)$.}

Let ${\sigma} \in {\Sigma}(C)$. The map ${\sigma}$ defines
 the effective
divisor
$L({\sigma}) = {\sigma}(x_{1}) + ... + {\sigma}(x_{18})$.
Clearly
$L({\sigma}) \in p_{*}^{-1}(C) = Supp \ {\Theta} \cup Supp \ P^{-}.$

\newpage
{\bf (2.2.4).}
{\it The sets ${\Sigma}_{\Theta}(C)$ and ${\Sigma}_{P^{-}}(C)$.}

\smallskip
{\bf Lemma.}
{\sl
Let $C$ be as above. Then

(1). The preimage $p_{*}^{-1}(C)$ coincides with the union of the
disjoint sets (each -- of cardinality $2^{17}$:

${\Sigma}_{\Theta}(C) =
\{ {\sigma} \in {\Sigma}(C) : L({\sigma}) \in Supp \ {\Theta} \}$
and

${\Sigma}_{P^{-}}(C) =
\{ {\sigma} \in  {\Sigma}(C) : L({\sigma}) \in Supp \ P^{-} \}.$

(2). $sgn({\sigma}', {\sigma}'') = +1$ iff both ${\sigma}'$ and
${\sigma}''$ belong to one of these two sets; otherwise $sgn = -1.$
}
\smallskip

{\bf Proof.}
Let $L({\sigma}')$ and $L({\sigma}'')$ be two elements of
$p_{*}^{-1}(C)$. The divisors $L({\sigma}')$ and $L({\sigma}'')$
are obtained from each other by a finite number of replacements
of the type: $L \mapsto L + l - \overline{l}$, where
$l + \overline{l} = p^{-1}(x)$ for some $x \in {\Delta}$. In this
case $x \in \{ x_{1},..., x_{18} \}$, and $L({\sigma}')$ and
$L({\sigma}'')$ are effective. Moreover, $L({\sigma}')$ and
$L({\sigma}'')$ can be regarded as  {\it general}  elements of
$Supp \ {\Theta} \cup Supp \ P^{-}$ (the cubic $C$ is general).
Therefore, $h^{0}(L({\sigma}'))$ and $h^{0}(L({\sigma}''))$ can be
only 1 or 2 (see e.g. [W]). Now, the lemma is a direct consequence
of the following statement:

\smallskip
{\bf (2.2.4)(*).}
{\sl
Let $\tilde{\Delta}$ be a smooth curve with an involution
$l \leftrightarrow \overline{l}$ without fixed points.
 Let $L$ be an
effective (see (2.1)) invertible sheaf on $\tilde{\Delta}$,
 and let
$l \in \tilde{\Delta}$. Then

$h^{0}(L)$ - $h^{0}(L + l - \overline{l}) \in \{ +1,-1 \}$
}
(see e.g. [Sh, 3.14], where (*) has been proved under more general
conditions).
{\bf q.e.d.}

\smallskip
{\bf (2.2.5).}
Let
$U = \{ C - \ \mbox{a smooth plane cubic}:\  C \cap {\Delta} = \ 18
 \ \mbox{disjoint points} \}$, let $C \in U$, and let
$Supp \ {\Theta}(C) =
 \{ L({\sigma}): {\sigma} \in {\Sigma}_{\Theta}(C) \}.$
Let
$Supp \ {\Theta}^{U} =
 \bigcup \{ Supp \ {\Theta}(C): C \in U \}.$
Then the algebraic set
$Supp \ {\Theta} \in Symm^{18}(\tilde{\Delta})$ is,
in an obvious way,
a closure of the open subset $Supp \ {\Theta}^{U}$.
The same (up-to replacing of the notation) is true also for
$Supp(P^{-})$ and $Supp(P^{-})^{U}$.

\smallskip
{\bf (2.3).}
{\it The global invariants $e({\Theta})$ and $e(P^{-})$.}
\smallskip

{\bf (2.3.1).}
{\it The maps $\tilde{\sigma}: S(C) \rightarrow S(L({\sigma})).$}

Let $C \in U$ be as in (2.2.5), and let ${\sigma} \in {\Sigma}(C)$
(see
(2.2)). Up-to change of the involutive line ($\leftrightarrow$
 a point
of $\tilde{\Delta}$), we may assume that
$L({\sigma}) = l_{1} + ... + l_{18}$. The lines $l_{i}$
and $\overline{l}_{i}$ are (-1)-curves
on the non-minimal ruled surface $S(C) = p^{-1}(C)$. Therefore,
 ${\sigma}$ defines,
in a unique way, a morphism
$\tilde{\sigma}: S(C) \longrightarrow S(L({\sigma}))$,
where $\tilde{\sigma}$ is the blow-down of the 18-tuple
$\{ \overline{l}_{1}, ..., \overline{l}_{18} \}.$

\smallskip
{\bf (2.3.2).}
{\it The invariants $e({\Theta})$ and $e(P^{-})$. }

Let $e(L({\sigma}))$ be the invariant of the ruled surface
$p({\sigma}): S(L({\sigma})) \rightarrow C$ (see (1.3.2)).
This way, we define a map:

$\tilde{e}: Supp \ {\Theta}^{U} \cup Supp(P^{-})^{U}
 \longrightarrow {\bf
Z}$,
$\tilde{e}: L({\sigma}) \mapsto  e(L({\sigma}))$,
(see (2.2.5)).

It is standard that the map $\tilde{e}$ must take a constant
 value
on some open
subset of each of the components of its domain.
Therefore, there exists a pair $e({\Theta}) , e(P^{-})$,
 and a pair of
(possibly smaller) open subsets
$Supp \ {\Theta}^{op}$ and $Supp(P^{-})^{op}$, such that
$\tilde{e}(L) = e({\Theta})$  for any $L \in Supp \ {\Theta}^{op}$,
$\tilde{e}(L) = e(P^{-})$  for any $L \in Supp \ (P{-})^{op}.$

We shall find these two numbers.

\smallskip
{\bf (2.4).}
$e({\Theta}) = 0, e(P^{-}) = -1.$
\smallskip

{\bf (2.4.1)}
{\bf Lemma.}
{\sl
Let $e = e({\Theta})$ and $\overline{e} = e(P^{-})$ be as in (2.3.2).
Then $\mid e - \overline{e} \mid  = 1.$
}
\smallskip

{\bf Proof.}
Let $C \in U$ ( see (2.2.5)) be such that
$L({\sigma}) \in Supp \ {\Theta}^{op} \cup Supp(P^{-})^{op}$,
for any ${\sigma} \in {\Sigma}(C)$ ( see (2.2),(2.3)).
In particular,
$e(S(L({\sigma}))) = e$  for any
$L({\sigma}) \in Supp \ {\Theta}(C)$, and
$e(S(L({\sigma}))) = \overline{e}$  for any
$\ L({\sigma}) \in Supp \ P^{-}(C)$.
Let $x_{i}, l_{i}, \overline{l}_{i}$, etc. be as in (2.2),(2.3),
and let ${\sigma}'$ and ${\sigma}''$ be such that
${\sigma}'(x_{i}) = {\sigma}''(x_{i})$ for any $i$, except
$i = j, \  j - fixed$. Then the minimal models
$S' = S(L({\sigma}'))$ and $S'' = S(L({\sigma}''))$
are obtained from each other by a single elementary transformation
{\sl elm}, centered in the point
$z_{j} = l_{j} \cap \overline{l}_{j} \in  p^{-1}(x_{j}).$
Now, (2.4.1) follows from the following:

\smallskip
{\bf (*).}
{\bf Sublemma} (see [LN, Lemma 4.3],or [Se, Lemma 7]).
Let $S' \rightarrow C$ and $S'' \rightarrow C$ be two ruled
surfaces  over the smooth base curve $C$, and let
$S'' = elm_{P}(S'')$, where $elm_{P}$ denotes the elementary
transformation of $S'$ -- centered in the point $P \in S'$. Then:

(i).  If no minimal section of $S'$ (see (1.3.4)) passes through $P$,
      then $e(S'') = e(S') - 1$.

(ii). If a minimal section of $S'$ passes through $P$, then
      $e(S'') = e(S') + 1$.

{\bf q.e.d.}

\smallskip
{\bf (2.4.2)}
{\bf Lemma.}
$\overline{e} = -1$.
\smallskip

{\bf Proof.}
Assume the contrary, i.e. $\overline{e} \ge 0$.

Let ${\cal L} \in {\Theta}$ be general, and let
$\mid {\cal L} \mid = \{ L(t), t \in {\bf P}^{1} \}$
be the linear system of ${\cal L}$ (see (1.2.2)).
Just as in (2.1.3), the pencil $\{ L(t) \}$ defines
the rational pencil of plane cubics $\{ C(t) = C(L(t)) \}$.
Since ${\cal L} \in {\Theta}$ is general, the general curve
$C(t)$  of the pencil $\{ C(t) \}$  is a smooth plane cubic,
and the only degenerations of $\{ C(t) \}$ are a finite
number of nodal plane cubics. Moreover, there is a finite number
of plane cubics $C(t) \in \{ C(t) \}$, which are simply tangent to
${\Delta}$. Anyway, there is an open subset
$V_{\cal L} \subset {\bf P}^{1}$ such that
$L({\sigma}) \in (Supp \ {\Theta})^{op} \cup Supp(P^{-})^{op}$
for any  $t \in V_{\cal L}$  and for any
${\sigma} \in {\Sigma}(C(t))$ ( see (2.2.2) -- (2.2.4)).

Let $L(t) = l_{1}(t) + ... + l_{18}(t), t \in V_{\cal L}$,
and let $L_{j}(t) = L(t) + \overline{l}_{j}(t) - l_{j}(t),\
j = 1,...,18$.

Let $\overline{e}$ be different from $-1$. Then
$\overline{e} \ge 0$,
and any of the ruled surfaces
$S(L_{j}(t))$  has only a finite number of minimal sections
 (see (1.3.4)).
The case (1.3.4)(2.a) can be excluded, because of considerations
 involving
the supposed general position. We shall prove that:

\smallskip
{\bf (*).}
{\bf Sublemma.}
{\sl
If the bidegree (2,2) divisor $T$ is general, then
the general $S(L), L \in Supp \ {\Theta} \cup Supp \ P^{-}$,
 cannot be
of type (1.3.4)(2.a) or of type (1.3.4)(3).
}
\smallskip

{\bf Proof.} On the one hand, since the general plane cubic $C(t)$
is a general plane cubic, it is in a general position w.r. to
the discriminant sextic ${\Delta}$. The various minimal
 ruled models
of the surface $S(C)$ are obtained from each other by elementary
transformations, related to the 18-tuple of degenerated fibers
$\{ l_{i} + \overline{l}_{i} = p^{-1}(x_{i}), i = 1,...,18 \}$.
On the other hand, we can fix for a moment the smooth plane
cubic $C$. Since the general plane sextic can be represented as
a discriminant of (one of the conic bundle structures of) some
bidegree (2,2) divisor $T$ (see [Ve]), we can choose the 18-tuple
$\{ x_{1},..., x_{18} \}$ without any closed restrictions.
The rest repeats the proof of [Se, Lemma 12]. {\bf q.e.d.}
\smallskip

So, the surface $S(L_{j}(t))$ has a finite number (one or two)
 of minimal
sections (By assumption,the case $\overline{e} = -1$ has been
 excluded).
In particular, this assumption, together with (*), imply that
 the general
$S(L), L \in Supp \ P^{-}$ has one (or two) minimal sections.
After taking the closure, one can define $\overline{\cal C}$
 to be the
family of all these minimal sections (related to the component
 $P^{-}$).
The last and (*) imply (under the assumption
 $\overline{e} \ge 0$ )  that
the family $\overline{\cal C}$ is 9-dimensional. The general
curve $\overline{C} \in \overline{\cal C}$ is mapped,
via $p = p_{1}$,
isomorphically onto a smooth plane cubic. Therefore,
$deg( \overline{C}) = (3,d)$, where $deg$ is the bidegree map.
The straightforward check, based on the normal bundle
 sequence for the triple
$\overline{C} \subset p^{-1}(p(\overline{C})) \subset T$, imply
that the total degree of
$\overline{C} \in$ (the 9-dimensional family)$\overline{\cal C}$
can be 8 , or 9 (see also [ZR]). Therefore, d = 5, or 6.

Let  ${\overline{\cal C}}_{\cal L}$  be the completion
(in $\overline{\cal C}$)
of the family of minimal sections of the surfaces $S(L_{j}(t))$
(see above). By definition, the base of this family is an algebraic
curve $B$. If the  general $S(L), L \in Supp \ P^{-}$ has two
 (resp one)
minimal sections, the base $B = B_{\cal L}$ is, in a natural way,
 a 36-sheeted
(resp. -- a 18-sheeted) covering of the projective line.
Let $\overline{S}_{\cal L}$ be the union of the curves
$C \in {\overline{\cal C}}_{\cal L}$.
Clearly, $\overline{S}_{\cal L}$ is an effective divisor on $T$
 (the curves,
which
sweep $\overline{S}_{\cal L}$ out, form an 1-dimensional
 algebraic family
parameterized by the algebraic curve $B$).

The surface $\overline{S}_{\cal L}$ can be reducible, or not.
 The irreducible
components of this surface correspond to the
irreducible components of the base $B = B_{\cal L}$. Let $B_{o}$
 be one of
these irreducible components, let
${\overline{\cal C}}_{o} \rightarrow B_{o}$
be the corresponding irreducible family, and let
 $\overline{S}_{o}$ be the
corresponding
irreducible component of $\overline{S}_{\cal L}$.

The surface  $\overline{S}_{o}$  represents the element

$cl(\overline{S}_{o}) \in {\bf Pic}(T) = {\bf Z}.l + {\bf Z}.h$,
 where
$l = cl(p_{1}^{*}{\cal O}_{{\bf P}^{2}}(1)),$
$h = cl(p_{2}^{*}{\cal O}_{{\bf P}^{2}}(1)),$   $ p_{1} = p.$

Therefore
$cl(\overline{S}_{o}) = al + bh$  for some integers  $a, b.$

\smallskip
{\bf (**).}
{\bf Sublemma.}
{\sl
The integers $a$ and $b$ are non-negative.
}
\smallskip

{\bf Proof.}
Let, for example, $b \le 0$. Let $f$ be the general fiber of
$p = p_{1} : T \rightarrow {\bf P}^{2}$. Since
$\overline{S}_{o} \in \mid al + bh \mid$ is effective
and $f$ is not a fixed curve on $T$, the intersection
number $(f.{\overline{S}_{o})_{T}}$ must be non-negative.
 Therefore
$0 \le  (f.(al + bh))_{T} = a(f.l)_{T} + b(f.h)_{T} =
 b(l^{2}.h)_{T} = 2b \le 0$ -- contradiction. {\bf q.e.d.}
\smallskip

Let $\overline{\cal C}_{o} =
    \{ \overline{C}({\xi}): {\xi} \in B_{o} \}$. The curves
$\overline{C}({\xi})$ belong to the same algebraic family.
In particular, the integer
$k = ( \overline{C}({\xi}_{1}).
\overline{C}({\xi}_{2}))_{\overline{S}_{o}}
 = ( \overline{C}({\xi})^{2})_{\overline{S}_{o}}$ is a
 (non-negative)
constant, which does not depend on the particular choice of
${\xi}_{1}, {\xi}_{2}, {\xi}.$
The general curve $\overline{C}({\xi}) \subset \overline{S}_{o}$
is a smooth elliptic curve of bidegree  $(3,d), d \in \{ 5,6 \}$.
The formal adjunction takes place on the divisor
$\overline{S}_{o}$ on $T$. In particular,
$0 = deg(K_{\overline{C}({\xi})}) =
 (K_{\overline{S}_{o}} + \overline{C}({\xi}).{\overline{C}({\xi})}) =
 ((a-1)l + (b-1)h).{\overline{C}({\xi})} + k =
 3(a-1) + d(b-1) + k$,
where $a, b, k \ge 0$ and $d \in \{ 5,6 \}$.

{\sl Case 1. $d = 5.$}
Then $3a + 5b + k = 8 \Rightarrow a = b = 1, k = 0.$

{\sl Case 2. $d = 6.$}
Then $3a + 6b + k =9 \Rightarrow a = b = 1, k = 0.$

In both cases $cl( \overline{S}_{o}) = l + h$, i.e.
$\overline{S}_{o}$ is a hyperplane section of $T$.
Since $\overline{S}_{o}$ is irreducible, it is a
surface of type $K3$ on $T$. ( In fact, the general
choice of the initial sheaf ${\cal L} \in {\Theta}$,
and the check of the parameters, imply that
$\overline{S}_{o}$ has to be a general hyperplane
section of $T$.) The family
$\overline{\cal C}_{o} \rightarrow B_{o}$
is an irreducible family of curves on $\overline{S}_{o}$,
and the general member of the family is an elliptic curve.
Therefore, $\overline{\cal C}_{o}$ is a rational pencil, i.e.
$B_{o} \cong {\bf P}^{1}$.

Let $\overline{C}({\xi})$ be general. Then the plane cubic
$C({\xi}) = p(\overline{C}({\xi}))$ belongs to the set $U$
(see (2.2.5)). In particular, $\overline{C}({\xi})$ is a
minimal section of some minimal model of the surface
$S(C({\xi})) = p^{-1}(C({\xi}))$ (see (2.3.1)).(In fact,
the corresponding minimal model must have two minimal sections --
 see e.g. (*).) Let
$\overline{L}({\xi})$ be the corresponding effective
divisor on $\tilde{\Delta}$ (ibid.) Clearly,
$\overline{L}({\xi}) \in Supp \ P^{-}$. By taking the completion,
we obtain the rational family
$\{ \overline{L}({\xi}): {\xi} \in {\bf P}^{1} \}
 \subset Supp \ P^{-}.$
In particular, all the $\overline{L}({\xi})$ belong to a same
(nontrivial) linear system. As it follows from the definition,
this system must be the linear system of some of the
divisors $L_{j}(t)$ (see the beginning of the proof). Clearly,
 $L_{j}(t)$
represents the general element of $Supp(P^{-})$ ( the initial sheaf
${\cal L}$ is a general element of ${\Theta}$). This is a
 contradiction,
since the linear system of the general $L \in Supp \ P^{-}$
 is trivial
 (see [W]). Therefore, $\overline{e} = e(P^{-}) = -1$. Lemma (2.4.2)
is proved.

\smallskip
{\bf (2.4.3)}
{\bf Corollary.}
{\sl
$e = e({\Theta}) = 0, \  \overline{e} = e(P^{-}) = -1$.
}

\bigskip
\centerline{\sc 3.
The family of minimal sections ${\cal C}_{\theta}$.}
\smallskip

{\bf (3.1).}
{\sl Definition of ${\cal C}_{\theta}$. }
\smallskip

Let $L \in Supp \ {\Theta}$ be general, and let
$S(L)$ be the corresponding minimal model of $S(C(L))$ (see (2.3)).
It follows from (2.4.1)-(2.4.3) that $e(S(L)) = e = 0$, and
$S(L)$ is decomposable (see (2.4.2)(*), and (1.3.4)(2)). Let
$C^{+} = C^{+}(L)$ and $C^{-} = C^{-}(L)$ be the minimal sections of
$S(L)$. By definition,

${\cal C}_{\theta}
 = \ \mbox{(the closure of)} \{ C^{+}(L), C^{-}(L) :
  L \in Supp \ {\Theta} \ \mbox{is general} \}$,
where the term ``general'' can be defined in an obvious way.

Clearly, $dim({\cal C}_{\theta}) = 9$ (the family ${\cal C}_{\theta}$
is generically a finite covering of the 9-dimensional projective
space of plane cubics).
Since $p = p_{1}$ maps the general $C^{+} \in {\cal C}_{\theta}$
isomorphically onto a plane cubic, $deg(C^{+}) = (3,d)$ for
some integer $d$. We shall prove the following

\smallskip
{\bf (3.1.2)}
{\bf Proposition.}
{\sl Let ${\cal C}_{\theta}$ be the family of minimal sections
defined in
(3.1.1), and let $C^{+}$ be a general element of ${\cal C}_{\theta}.$
Then:

(a). $ deg(C^{+}) = (3,6)$;

(b). $C^{+}$ lies on a hyperplane section $S_{C^{+}} \subset T$.
}

\smallskip
{\bf Proof.}
Let $L = L(0) \in Supp \ {\Theta}$ be general, and let
$C^{+}(L) = C^{+}(0)$ and $C^{-}(L) = C^{-}(0)$ be the
minimal sections of $S(L) = S(L(0))$. Just as in the proof of
Lemma (2.4.2), we can prove that: {\bf (1).} The sections
 $C^{+}(0)$ and
$C^{-}(0)$ can be included in rational families
$ \{ C^{+}(t) : t \in {\bf P}^{1} \}$  and
$ \{ C^{-}(t) : t \in {\bf P}^{1} \}$.
{\bf (2).} The enveloping surfaces $S^{+} = S_{C^{+}}$ and
$S^{-} = S_{C^{-}}$ of these two families are hyperplane sections
of the threefold $T$.This proves (b).

\smallskip
{\bf (*).}
{\bf Remark.} The check of the parameters implies that
the general hyperplane section of $T$ carries such a pencil
of minimal sections.Moreover, it is not hard to see that
the general minimal section determines its enveloping hyperplane
sections in a unique way -- there is a unique hyperplane in
${\bf P}^{8}$  which contains the general curve
$C \in {\cal C}_{\theta}$
(see also [T] and [C], where corresponding
results are proved for the Reye sextics on the Quartic double
solid).

The part (a) is a corollary from (b), and from standard
arguments involving
the
normal bundle sequences for the triples
$C^{+} \subset S(C) \subset T$ and
$C^{+} \subset S^{+} \subset T$ (similarly -- for $C^{-}$).
 {\bf q.e.d.}

\bigskip
\centerline{\sc 4.
The Abel-Jacobi image of the family ${\cal C}_{\theta}$.}
\smallskip

{\bf (4.1)}
{\bf Theorem.}
{\sl
Let ${\Phi} : {\cal C}_{\theta} \longrightarrow J(T)$
be the Abel-Jacobi map for the family ${\cal C}_{\theta}$
}
(see e.g. [C]).
{\sl
Then the image ${\Phi}({\cal C}_{\theta})$
is a copy of the theta divisor ${\Theta}(T)$.
}
\smallskip

{\bf Proof.}
It follows from (3.1.1),(3.1.2) that the family
${\cal C}_{\theta}$ is generically a 2-sheeted
covering of the family $Supp \ {\Theta}$ :
(Let $L \in Supp \ {\Theta}$ be general. Then the
fiber ${\xi}^{-1}(L)$ of the covering
${\xi} : {\cal C}_{\theta} \rightarrow Supp \ {\Theta}$
is actually the pair $C^{+}(L),C^{-}(L)$ -- see above).
It follows from the definition of $Supp \ {\Theta}$ (see (2.1))
that the (Prym-)Abel-Jacobi map
${\Psi}$ sends $Supp \ {\Theta}$ onto the copy ${\Theta}$
described in $(1.2.2)$. In fact, ${\Psi}$ contracts
the linear systems $\mid {\cal L} \mid$, ${\cal L} \in {\Theta}$,
to the points ${\cal L}$.

Let $L = L(0) \in Supp \ {\Theta}$ be general, and let
$C^{+}(L) = C^{+}(0)$ and $C^{-}(L) = C^{-}(0)$ be the minimal
sections of $S(L) = S(0)$. Let
$\{ C^{+}(t) \}$ and $\{ C^{-}(t) \}$ be the rational families
defined as in the proof of (3.1.2), and let
$\{ L^{+}(t) \}$ and $\{ L^{-}(t) \}$ be the corresponding
pencils in $Supp \ {\Theta}$. Since $L^{+}(0) = L^{-}(0) = L(0) = L$,
these two pencils coincide. Clearly,
$L^{+}(t) = L^{-}(t) = L(t)$, where $\{ L(t) \}$ is the pencil
defined by the linear system
$\mid {\cal L} \mid  = \mid  {\cal O}_{\tilde{\Delta}}(L) \mid$.
The general $C^{+}(t)$ and $C^{-}(t)$ are the minimal sections
of the corresponding surfaces $S(L(t))$.
Therefore, the natural 2-sheeted covering ${\xi}$ factors through
general linear systems $\mid {\cal L} \mid , {\cal L} \in {\Theta}$.
Now, the only which is left, is to prove the following

\smallskip
{\bf (4.2)}
{\bf Lemma.}
{\sl
Let  $L \in Supp \ {\Theta}$ be general, and let
$C^{+} = C^{+}(L)$ and
$C^{-} = C^{-}(L)$ be the minimal sections of $S(L)$. Then
${\Phi}(C^{+}(L)) = {\Phi}(C^{-}(L))$.
}
\smallskip

{\bf Proof} of the Lemma.
Let $C = p(C^{+}) = p(C^{-})$ be the isomorphic image of
$C^{+}$ and $C^{-}$, and let $S(C) = p^{-1}(C) \subset T$.
Let $S(L) = {\bf P}_{C}({\cal O} \oplus \varepsilon )$ be as
in (1.3.4.)(2). By definition, $C^{+}$ is the tautological
section of ${\cal E}^{+} = {\cal O} \oplus \varepsilon$, and
$C^{-}$ is the tautological section of
${\cal E}^{-} = {\cal O} \oplus {- \varepsilon}$.
Let $p_{L} : S(L) \rightarrow C$ be the projection, and let
$\sim$ denotes the linear equivalence on the ruled surface
$S(L)$. It follows from the definitions of $C^{+}$ and
$C^{-}$ that
$C^{+} \sim C^{-} + {\varepsilon}.f_{L}$
where ${\varepsilon}.f_{L} := p_{L}^{*}({\varepsilon})$.

Let ${\delta}(C) = {\Delta} \cap C$, and let $Q$ be a point
of $C$. Since $C$ is an elliptic
curve and ${\varepsilon} \in {\bf Pic}^{o}(C)$, the sheaf
${\varepsilon}$ can be represented in the form
${\varepsilon} = {\cal O}_{\Delta}(P - Q)$,
where $P = P({\varepsilon})$ is determined uniquely by the
choice of $Q$ (the sheaf ${\varepsilon}$ is fixed).
Clearly, because of the continuous choice of the initial
point $Q$, we can always choose $Q$ such that
$\{ Q,P \} \cap {\delta}(C) = \oslash$. Then
$C^{+} \sim C^{-} + (P - Q).f_{L}$, and $f_{L}$ can be
also regarded as a fiber of the surface
$S(C) \subset T$  (equivalently -- $f_{L}$ can be replaced
by the general fiber $f$ of the projection
$p : T \rightarrow {\bf P}^{2}$ ).
Therefore ${\Phi}(C^{+}) = {\Phi}(C^{-}  + (P - Q).f) =
 {\Phi}(C^{-})$, since the fibers $P.f = p^{-1}(P)$
and $Q.f = p^{-1}(Q)$ are rationally equivalent on $T$.
{\bf q.e.d.}

The Theorem is proved.

\bigskip
\centerline{\sc 5.
Hyperplane sections of the Verra threefold,}

\centerline{\sc and plane
sextics with vanishing theta-null.}
\smallskip

{\bf (5.1)}
{\bf Proposition.}
{\sl
Let $T$ be a smooth bidegree (2,2) divisor in
${\bf P}^{2} \times {\bf P}^{2}$
(a Verra threefold), let $T \subset {\bf P}^{8}$
be the  Segre embedding, and let
$p_{i}:T \rightarrow {\bf P}^{2}, i = 1,2 \ (p_{1} = p)$
be the canonical projections. Let $S \subset T$ be
a general (esp. -- irreducible) hyperplane section of $T$.
Then:

(1). $p_{i}:S \rightarrow {\bf P}^{2}$ is a double covering
branched along a plane sextic $D_{i}$ which is totally tangent
to the discriminant sextic ${\Delta}_{i}$, i = 1,2.

(2). $D_{i}$ has a vanishing theta-null, i = 1,2.
}
\smallskip

{\bf Proof.}
On the one hand, the family of hyperplane sections of $T$
 is naturally
isomorphic to $({\bf P}^{8})^{*}$, and the general hyperplane
 section
of $T$ is a smooth $K3$ surface. On the other hand,
 the Abel-Jacobi
map sends the 9-dimensional family ${\cal C}_{\theta}$ onto a copy
of the 8-dimensional theta divisor ${\Theta} = {\Theta}(T)$,
and the general $C \in {\cal C}_{\theta}$ lies in a unique
hyperplane
section. On the other hand, the components of the general fiber of
${\Phi}:{\cal C}_{\theta} \rightarrow {\Theta}$ are (pairs of)
 pencils
of minimal sections, and the general pencil of this type is actually
an elliptic pencil on some hyperplane section  $S$ of $T$
(see (3.1.2)(*) and Theorem 4.1).
A comparison between the dimensions, and standard arguments involving
general positions, imply that the general hyperplane section of $T$
(which is a smooth $K3$ surface)
carries a finite number of pencils of minimal sections.
Obviously, this argument can be applied both: to $p = p_{1}$, and
to $p_{2}$.
Clearly, $p_{i}:S \rightarrow {\bf P}^{2}$ is a double covering.
Since $S$ is a $K3$ surface, the branch locus $D_{i}$ of $p_{i}$
is a sextic curve, and it is not hard to see that $D_{1}$
is totally tangent to the discriminant sextic ${\Delta}_{i}$
$( \Rightarrow (1).).$

Fix for a moment $p_{i} = p_{1}$.
Let $\{ C(t) : t \in {\bf P}^{1} \}$ be one of the
elliptic pencils of minimal sections (w.r. to $p_{1}$) on $S$.
The elements of this pencil are curves of arithmetical genus 1
and of bidegree (3,6), and  $p_{1}$  projects the pencil
$\{ C(t) \}$ onto a pencil $\{ p_{1}(C(t)) \}$ of plane cubics.
Clearly, $p_{1}(C(t))$ is totally tangent to the branch locus
$D_{1}$ of the double covering $p_{1}:S \rightarrow {\bf P}^{2}$.
Therefore, $\{ p_{1}(C(t)) \}$ defines a vanishing theta-null
${\chi}_{1}$ on $D_{1}$. Similarly -- for $p_{2}$.
$( \Rightarrow (2).).$

\smallskip
{\bf (5.2)}
{\bf Remarks.}
\smallskip

{\bf (5.2.1)}
{\it
Prym-Canonical systems of $({\Delta}_{i},{\eta}_{i})$,
and sextics with vanishing theta-null.
}

Let $i \in \{ 1,2 \}$, and let ${\eta}_{i}$ be the torsion sheaf
which defines the double covering
$ \overline{\Delta}_{i} \rightarrow {\Delta}_{i}$.
Since  $D_{i}$ is totally tangent to ${\Delta}_{i}$,  \
$D_{i}\mid_{{\Delta}_{i}} = 2.{\delta}_{i}$, for some
effective divisor ${\delta}_{i} = {\delta}_{i}(S)$ on
${\Delta}_{i}$. It can be seen (see [Ve]) that
${\delta}_{i}$ belongs to the Prym-Canonical linear system
${\omega}_{{\Delta}_{i}}({\eta}_{i}) =
{\cal O}_{{\Delta}_{i}}(3) \otimes {\eta}_{i}$.
However, the effective divisor ${\delta}_{i}$ does not determine
the totally tangent sextic (along ${\delta}_{i}$) in a
unique way. In fact, any plane sextic ($ \not= {\Delta}_{i}$)
of the pencil $\{ D_{i}(t) \}$, which is  spanned on  $D_{i}$  and
${\Delta}_{i}$,
is totally tangent to ${\Delta}_{i}$ along ${\delta}_{i}(S)$.
Proposition (5.1) tells that the branch sextic $D_{i} = D_{i}(S)$
``chooses'', among the curves of this pencil, a sextic with a
vanishing theta-null.

\smallskip
{\bf (5.2.2)}
{\it
The two splittings.
}

In fact, the non-uniqueness of the minimal section of the
general surface $S(L), L \in Supp \ {\Theta}$ (see e.g. (3.1)),
implies that the family of plane cubics $\{ p_{1}(C(t)) \}$
(see the proof of (5.1)) defines a vanishing theta-null
on at least two branch loci: the branch loci of the
enveloping hyperplane sections $S^{+}$ and $S^{-}$
(see the proof of Proposition (3.1.2)).

A similar splitting takes place also on the general hyperplane
section
$S$. Let e.g. $\{ C(t) \} \}$ be one of the pencils of minimal
sections on the hyperplane section $S$ ; we let $p = p_{1}$.
Let $\{ p(C(t) \}$ be the corresponding rational pencil of
plane cubics, and let $\{ S(t) = p^{-1}(p(C(t)) \}$ be the pencil
of canonical conic bundle surfaces above $\{ p(C(t)) \}$.
Let $S$ be the enveloping hyperplane section of $\{ C(t) \}$.
Then the pencil
$\{ S \cap S(t) \}$
(of singular curves of arithmetical genus 10, and
of bidegree (6,12) )  splits into two pencils:
 the pencil $\{ C(t) \}$, and the  {\it conjugate}  pencil
$\{ \overline{C}(t) = S \cap S(t) - C(t) \}$.
Obviously, $\{ \overline{C}(t) \} \subset {\cal C}_{\theta}$,
and the two pencils $\{ C(t) \}$ and $\{ \overline{C}(t) \}$
(which lie on the same hyperplane section $S$)
are projected onto the same pencil of plane cubics.

\smallskip
{\bf (5.2.3)}
{\it
The Double Flag, and the Dixon correspondence between the
branched loci $D_{1}(S)$ and $D_{2}(S)$.
}

Let $T, {\Delta}_{i}, {\eta}_{i}$, etc., be as usual.
The identity between the Prym varieties
$P({\Delta}_{1},{\eta}_{1})$ and
$P({\Delta}_{2},{\eta}_{2})$ gives a counterexample to
the Torelli theorem for Prym varieties ([Ve]).
There exists another family of counterexamples to
the Torelli theorem for Prym varieties, which can be
regarded as a ``hyperelliptic'' degeneration of the
family of bidegree (2,2) divisors. This is the family of
Double Flags (see below),  and the general member of this family
can be described as follows:

Let $seg:{\bf P}^{2} \times {\bf P}^{2} \rightarrow {\bf P}^{8}$
be the Segre embedding, and let $W$ be the image of $seg$.
Then the general hyperplane section $Y$ of $W$
is a smooth Fano threefold, and any such a threefold
is isomorphic to the Flag variety ${\bf P}(T_{{\bf P}^{2}})$ --
the incidence correspondence between points and lines on the
projective plane ${\bf P}^{2}$.
By definition, a {\it Double Flag} is any double covering
$\pi:X \rightarrow Y$, branched along an intersection $S$ of $Y$
and a quadric. Clearly, the smooth {\it Double Flag} $X$
is a smooth
Fano threefold, and the maps
$\tilde{p}_{i} = p_{i} \circ {\pi}: X \rightarrow {\bf P}^{2}$
\ $(i = 1,2)$ define conic bundle structures on $X$. The
discriminant curve of $\tilde{p}_{i}$ coincides with the
branch locus $D_{i} = D_{i}(S)$ of the 2-sheeted covering
$p_{i}:S \rightarrow {\bf P}^{2}$, i = 1,2. Obviously,
the branch surface $S$ is also a hyperplane section of
a Verra threefold. In particular, the discriminant sextic
$D_{i}$ has a vanishing theta-null. Just as in the case of
the Verra threefold, the conic bundle
structure $\tilde{p}_{i}$ defines, in a standard way,
the torsion-2 sheaf ${\tilde{\eta}}_{i}$ on $D_{i}$.
In [Ve] A.Verra proves that the discriminant pairs
$({\Delta}_{1},{\eta}_{1})$ and
$({\Delta}_{2},{\eta}_{2})$,
for the bidegree (2,2) divisor $T$,
are connected to each other by the classical Dixon
correspondence. The same argument can be applied to
the couple of discriminant pairs for the Double Flag.
Remember that the discriminant sextics for the Double Flag
are also branched loci, related to a hyperplane section
of a Verra threefold. Obviously, the general
hyperplane section $S$ of $T$ is also a branch locus
of a double covering of the Flag variety $Y$.
Therefore, the branch loci $D_{1}(S)$ and $D_{2}(S)$
(defined by the general hyperplane section $S$ of the
threefold $T$)  can be included in discriminant pairs
which are connected by the Dixon correspondence.
In fact, the condition imposed by the existence of a
theta-null, is a closed condition of codimension 1,
on the 19-dimensional space of non-isomorphic
plane sextics. The K3-surfaces, which are double coverings
of ${\bf P}^{2}$, also form a 19-dimensional moduli space
${\cal M}$. A vanishing theta-null, imposed on the branch
locus of such a K3-surface,  separates a codim.1 subspace
${\cal M}_{o}$ of ${\cal M}$.  Note also that the Double
Flags are codim.1 degeneration of Verra threefolds, and the
general plane sextic appears as one of the discriminant curves
for some Verra threefold.  This implies the following:

\smallskip
{\bf (5.3)}
{\bf Corollary.}
{\sl
There exists a component (of maximal dimension)
${\cal M}'$ of the 18-dimensional parameter space ${\cal M}_{o}$
(of K3-surfaces which are double coverings of ${\bf P}^{2}$
branched along sextics with a vanishing theta-null), such that
the general element $S$ of ${\cal M}'$ can be represented
by two ways as a double covering of ${\bf P}^{2}$. Moreover,
if $D_{1}$ and $D_{2}$ are the branch loci of these two
coverings, then $D_{1}$ and $D_{2}$ can be included into
discriminant pairs related to each other by the Dixon
correspondence.
}

\bigskip
\centerline{\sc 6.
The family ${\cal D}$ of elliptic sextics of
bidegree (3,3) on $T(2,2)$.
}
\smallskip

{\bf (6.0).}
In this section we prove that the Abel-Jacobi map sends
the 6-dimensional family ${\cal D}$ of elliptic sextics
on $T$ of bidegree (3,3) onto a 3-dimensional subvariety
$Z \subset J(T)$ -- isomorphic to a component of
$Sing \ {\Theta}$.

\smallskip
{\bf (6.1).}
{\it
Definition of ${\cal D}$.
}

Let $C \subset T$ be a reduced and connected curve of
arithmetical genus $p_{a}(C) = 1$, and of total degree
$\mid deg(C) \mid  \ = 6$. Here $deg$ is the usual bidegree
map $deg:\ \{ \mbox{1-cycles on} \ T \} \rightarrow
{\bf Z} \oplus {\bf Z}$.
If such a curve $C \subset T$ does exist, the bidegree
$deg(C)$ can be one of the pairs (2,4), (3,3), (4,2).
We shall consider the middle case: $deg(C) = (3,3)$.
By definition:

${\cal D} = \ \mbox{(the closure of)}
\{ C - \ \mbox{a reduced and connected curve on} \ T :
 p_{a}(C) = 1, deg(C) = (3,3) \}$.

We call ${\cal D}$ the
{\it family of elliptic sextics of bidegree (3,3)}
on $T$.

\newpage
{\bf (6.2).}
{\it
The existence of $C \in {\cal D}$.
}
\smallskip

We shall find at least one reduced and connected curve
$C \in {\cal D}$.The existence of such a curve makes it
possible to apply general techniques (see e.g. [ZR])
to obtain more information about the non-empty family
${\cal D}$.

\smallskip
{\bf (6.2.1).}
{\it
Curves of small degree on the components of the reducible
hyperplane sections of $T$.
}

Let $l \subset {\bf P}^{2}$ be a {\it general} line, and let
$S_{l} = p^{-1}(l)$;  here, as usual,  $p = p_{1}$, etc.
The following can be seen immediately:

  {\bf (1).} $deg(S_{l}) = 6; K_{S_{l}} = -h$ ; here
$h$  is the restriction on  $S_{l}$  of the class
$h = cl(p_{2}^{*}{\cal O}_{{\bf P}^{2}}(1)).$

  {\bf (2).} Let $q = p_{2}$ be the  2-nd projection.
Then $q:S_{l} \rightarrow {\bf P}^{2}$ is a double
covering branched along a smooth quartic curve.
In particular, $S_{l}$ is isomorphic to a del Pezzo
surface with $K^{2} = 2$, and the projection $q$
defines the anticanonical linear system on $S_{l}$.

  {\bf (3).} Let ${\Delta} = {\Delta}_{1}$ be the
discriminant sextic for $p$, and let
${\Delta} \cap l = \{ {\xi}_{1},...,{\xi}_{6} \}$
( $l$ is general). The projection $p$ separates
one of the various conic bundle structures on
the del Pezzo surface $S_{l}$. In particular,
$p:S_{l} \rightarrow l$ defines (in a non-unique way)
a morphism ${\sigma}:S_{l} \rightarrow {\bf P}^{2}$, s.t.:

  (i). ${\sigma}$ is a composition of seven  ${\sigma}$-processes
${\sigma}_{x_{i}}, i=0,...,6.$  Here
$x_{0} ,..., x_{6}$  are 7 points in  ${\bf P}^{2}$, in a general
position w.r.to the plane curves of degree  $\le 3$
(see [DPT]).

  (ii). The pencil of
conics $p:S_{l} \rightarrow l$ is defined by the
${\sigma}$-preimage of the non-complete linear system
$\mid {\cal O}_{{\bf P}^{2}}(1 - x_{o}) \mid$.

  (iii). The degenerate fiber
$p^{-1}({\xi}_{i}) = [x_{i}] + [x_{o},x_{i}]$,
where $[x_{i}]$ is the exceptional curve of
${\sigma}_{i}$, and $[x_{o},x_{i}]$ is the proper
${\sigma}$-preimage of the line $<x_{o},x_{i}>$,
i=1,...,6.

  {\bf (4).} The (-1)-curves $[x_{i}]$ and $[x_{o},x_{i}]$,
i=1,...,6, are the only lines on the surface
$S_{l} \subset {\bf P}^{5}(l) = Span(S_{l}) \subset {\bf P}^{8}$,
and the embedding $S_{l} \subset {\bf P}^{5}(l)$ is defined
by the non-complete linear system
$\mid
{\cal O}_{{\bf P}^{2}}( 4 - 2x_{o} - (x_{1} + ... + x_{6}) )
\mid$.

The next lemma is a direct corollary of the properties
(1),(2),(3),(4). We refer to the articles I -- V of
M.Demazure in [DPT],
which provide a comprehensive study of the del Pezzo surfaces.

\smallskip
{\bf (*).Lemma.}
{\sl
Let $S_{l} = p^{-1}(l)$ be general. Then, in the terms of (1)-(4):

{\bf (a).}
There exist exactly $2^{5} = 32$ morphisms
$S_{l} \longrightarrow {\bf F}_{1} =
 {\bf P}_{{\bf P}^{1}}({\cal O} \oplus {\cal O}(1))$,
and exactly 32 morphisms
$S_{l} \longrightarrow {\bf F}_{o} =
 {\bf P}({\cal O} \oplus {\cal O}(1))$,
in which $l$ remains a base, and $l_{i} = [x_{i}]$ and
$\overline{l}_{i}$, i=1,...,6, are mapped to fibers or points.

{\bf (b).}
The 32 curves on $S_{l}$, which are mapped onto the
exceptional curves of the 32 ruled surfaces ${\bf F}_{1}$,
are conics of bidegree (1,1). These 32 conics are:

  (i). the (-1)-curve $[x_{o}]$;

  (ii). the (-1)-curves $[x_{o},x_{i}]$, i=1,...,6;

  (iii). the (-1)-curves
$ \left[ x_{o},...,\widehat{x_{i}},...,
\widehat{x_{j}},...,x_{6} \right] $  :=
(the proper preimage of the conic through the corresponding
5-tuple of points);

  (iv). the (-1)-curve
$[2x_{o}, x_{1},...,x_{6}]$ := (the proper preimage of the
cubic through $x_{1},...,x_{6}$, which has a node in $x_{o}$).

Moreover, these 32 curves are the only bidegree (1,1)-conics on
$S_{l}$.

{\bf (c).}
The general curve $C$ on $S_{l}$, which is mapped onto a
0-section on ${\bf F}_{o} \rightarrow {\bf P}^{1} = l$  (!),
is a twisted cubic on $S_{l}$ of bidegree (1,2). The 32
pencils of (1,2)-curves on $S_{l}$, which correspond to
the 32 fibrations of type ${\bf F}_{1}$, can be described
explicitly -- as in (b).  Just as in (b), these 32 pencils
of twisted cubics, are the only pencils of bidegree
(1,2)-curves on $S_{l}$, and any bidegree (1,2)-curve on
$S_{l}$ is an element of one of these pencils.
}

\smallskip
{\bf (**). Corollary.}
{\sl

  {\bf (a).}
The family ${\cal C}_{1,1}$ of conics on $T$ of bidegree (1,1)
is 2-dimensional. (In fact, ${\cal C}_{1,1}$ is, in a natural way,
a 32-sheeted covering of ${\bf P}^{2}$ -- see above). There are
finitely many conics on $T$ of bidegree (1,1), which pass
through the general point of $T$.

  {\bf (b).}
The family ${\cal C}_{1,2}$ of twisted cubics on $T$ of
bidegree (1,2) is 3-dimensional. If $z \in T$ is a general
point, then the family ${\cal C}_{1,2}(z)$ of these
bidegree (1,2)-cubics on $T$ which pass through $z$
is an algebraic set of dimension 1.
}

\smallskip
{\bf Proof.}
Standard.
\smallskip

{\bf (6.2.2).}
{\it
The existence of a reduced curve $C \in {\cal D}$.
}

We shall find such a curve.

Let $l \subset {\bf P}^{2}$ be a general line, and let
$C_{1,1} \subset T$ be one of the 32 conics of bidegree (1,1)
 ``above'' $l$ (see (6.1.1)(*)(b)). The 2-nd projection
$q = p_{2}$ maps $C_{1,1}$  isomorphically onto a line
$l' \subset {\bf P}^{2}$. The line $l'$ intersects the
2-nd discriminant sextic ${\Delta}_{2}$ in 6 points.
Therefore, there are 6 lines on $T$ of bidegree (1,0),
which intersect the curve $C_{1,1}$. Let $L_{1,0}$ be
one of these lines,and let $m = p(L_{1,0})$.
Let $x$ be a general point of $C_{1,1}$. Clearly,
$x$ can be regarded as a general point of  $T$.
Let $S(x)$ be the union of all the curves which belong
to the 1-dimensional family ${\cal C}_{1,2}(x)$
(see (6.2.1)(**)).  $S(x) \subset T$ is a surface,
and it is not hard to see that the general choice
of $x \in C_{1,1}$ implies that $S(x)$ does not
contain the line $L_{1,0}$. Therefore, $S(x)$ intersects
$L_{1,0}$ in a finite number of points (different form the
point $C_{1,1} \cap L_{1,0} .$  Let $y \in L_{1,0}$
be one of these points. In particular, there exists a
twisted cubic $C_{1,2} \in {\cal C}_{1,2}$, which passes
through $x$ and $y$. Obviously, the curve
$C = C_{1,1} + L_{1,0} + C_{1,2}$ is reduced, and
$deg(C) = (3,3)$.

The existence of $C \in {\cal D}$ assures the non-emptiness
of ${\cal D}$.  The general results in [ZR] imply that
the family ${\cal D}$ is at least  6-dimensional.  It is
not hard to see that the general curve $C \in {\cal D}$
is smooth, and the projections $p_{1}$ and $p_{2}$
send $C$ izomorphically onto the smooth plane cubics
$C_{1}$ and $C_{2}$. Let $S_{C_{1}} = p_{1}^{-1}(C_{1})$ and
$S_{C_{2}} = p_{2}^{-1}(C_{2})$. Just as in sections 2, 3, 4, the
curve $C$ defines, via intersection,  the effective divisors
$L_{i}(C)$ on $\tilde{\Delta}_{i}$,
which belong to the linear system of some
${\cal L}_{i} \in Nm^{-1}({\omega}_{{\Delta}_{i}}), \ i = 1,2$.

According to the agreement, we fix $p = p_{1}$.
Let ${\Delta} = {\Delta}_{1}, \ L(C) = L_{1}(C)$, etc.,
be the corresponding objects. Obviously, $C$ is a minimal
section of the ruled surface $S(L(C))$ (see (2.3.1)), and
the same arguments  as in sect. 2, 3, 4  imply that
$S(L(C))$ is indecomposable and $e(S(L(C))) \ge 3$. (The
same normal bundle sequence argument gives
$dim({\cal D}) = 6$).
In particular, $C$ is the unique minimal section
of the canonical surface $S(p(C)) = p^{-1}(p(C))$
(which corresponds to the
unique minimal section of the ruled model $S(L(C))$).

The curve $C$ is not an element of the total family
of minimal sections  ${\cal C}_{\theta}$.
This is caused by the special position of the cubics
$p(C) , C \in {\cal D}$,  which form a 6-dimensional
subset of the 9-space of all the plane cubics.
We collect the last in the following

\smallskip
{\bf (6.2.3)}
{\bf Lemma.}
{\sl

  {\bf (i).}  $dim{\cal D} = 6$;

  {\bf (ii).} The general curve $C \in {\cal D}$
is a smooth sextic
of bidegree (3,3), and the projections $p_{1} = p$ and $p_{2}$
map $C$ isomorphically onto plane cubics. In particular:

  {\bf (iii).} ($p = p_{1}$).
Let $L(C) \in Symm^{18}(\tilde{\Delta})$
be, as usual, the effective divisor on $\tilde{\Delta}$,
defined by the 18-tuple of lines of bidegree (0,1) which
intersect the curve $C$. Then
$L(C) \in Supp \ {\Theta} \cup Supp \ P^{-}$.
}

\newpage
{\bf (6.3).}
{\it
The family $\cal D$ and the quadrics of rank 6 through $T$.
}
\smallskip

{\bf (6.3.1).}
{\it
The elements of ${\cal D}$ as components of canonical
curves on $T$.
}

Let ${\bf P}^{6}$ be a subspace of ${\bf P}^{8} = Span(T)$
such that $dim(T \cap {\bf P}^{6}) = 1$, and let
$C({\bf P}^{6}) = T \cap {\bf P}^{6}$. The curve
$C({\bf P}^{6}) \subset T$ is a canonical curve of
degree 12 , and of arithmetical genus 7.  Call such a curve $C$
a {\it canonical curve on $T$ }.  Obviously, all the
canonical curves on $T$ are rationally equivalent,
and the family of canonical curves on $T$
can be represented by the 14-dimensional Grassmann variety
$G(7,9) = G(6:{\bf P}^{8})$.

The the curves of the family ${\cal D}$ are closely connected
with the degenerations of the family of canonical curves on
$T$.  More precisely, let $C \in {\cal D}$ be general, let
${\bf P}^{5}(C) = Span(C)$, and let
${\bf P}^{6} \in {\bf P}^{2}(C) := {\bf P}^{8}/{\bf P}^{5}(C)$
be a 6-space through ${\bf P}^{5}(C)$.  Then the canonical
curve $C({\bf P}^{6})$ splits into two components:
 $C({\bf P}^{6}) = C + \tilde{C}$
where $\tilde{C} \in {\cal D}$, and
${\delta}(C,\tilde{C}) = {\#}(C \cap \tilde{C}) = 6$.

\smallskip
{\bf (6.3.2).}
{\it
The determinantal subvarieties of  ${\bf I}_{2}(W).$
}

Let $W \subset {\bf P}^{8}$ be the Segre image
of  ${\bf P}^{2} \times {\bf P}^{2}$, and let
${\bf P}^{2} \times {\bf P}^{2}
 = {\bf P}(E) \times {\bf P}(F)$, where
$E$ and $F$ are complex 3-spaces.
Let ${\bf P}^{8} = {\bf P}(E \otimes F)$.
Then, as it follows from the definition of the Segre map,
the elements of $W \subset {\bf P}^{8}$ are in 1:1 correspondence
with the ${\bf C}^{*}$-classes of unitary tensor products
$u \otimes v : u \in E, v \in F$. In particular, let
$(x_{i},e_{i})$ and $(y_{j},f_{j})$
be any coordinate systems on $E$ and $F$.
Then  $(z_{ij} = x_{i}.y_{j} \ , \ g_{ij} = e_{i} \otimes f_{j})$
is
a coordinate system on $E \otimes F$. Let $[z_{ij}]$ be the
coordinate matrix, and let ${\bf I}_{2}(W)$ be the projective
space of quadrics in ${\bf P}^{8}$ which pass through $W$.
Then, in matrix coordinates $z_{ij}$,
${\bf I}_{2}(W)$
is spanned on the 9 quadratic equations
$rank([z_{ij}]) = 1$, i.e. ${\bf I}_{2}(W)$ is a projective
8-space. Clearly, the choice of the coordinates $z_{ij}$
defines a linear isomorphism

${\psi}(z):{\bf P}^{8} = {\bf P}(E \otimes F) \cong {\bf I}_{2}(W)$.

The linear map  ${\psi}(z)$  sends the ${\bf C}^{*}$-classes
of the unit tensor products (i.e. -- the elements of $W$)
 to quadrics of rank 4.  Moreover, any quadric of rank 4,
which contains $W$, can be represented in this way.
 It is well-known (see e.g. [LV]) that the set

$Sec(W)$ := the closure of the union of all the bisecant
lines of $W$,

is a cubic hypersurface in ${\bf P}^{8}$,  i.e., $W$
is one of the four
Severi varieties (ibid.). It follows from the definition
of ${\psi}(z)$
that ${\psi}(z)$ sends the points of $Sec(W)$ to quadrics of
rank 6.  It is not hard to see that any quadric of rank 6
which contains $W$  can be represented (in a non-unique way)
as an image of a point of $Sec(W)$. We collect these observations
in the following:

\smallskip
{\bf (*).
Lemma}.
{\sl
Let ${\bf I}_{2}(W)$ be the projective space of quadrics in
${\bf P}^{8}$ which contain the fourfold $W$, and let
$D_{k}(W)  = D_{k}({\bf I}_{2}(W)) =
\ \mbox{(the closure of)} \{ Q \in {\bf I}_{2}(W) : rank(Q) =  k \}$
be the $k$-th determinantal of ${\bf I}_{2}(W)$.
Then $D_{k}(W) \neq \oslash  \Leftrightarrow  k \in \{ 4,6,9 \}$.
Moreover, there exists a linear isomorphism
${\psi}: {\bf P}^{8} \rightarrow {\bf I}_{2}(W)$
such that

${\psi}(W) = D_{4}(W)$,

${\psi}(Sec(W)) = D_{6}(W).$
}
\smallskip

{\bf Proof.}
Let $[z_{ij}]$ be as above, and let ${\psi} = {\psi}(z)$.
Then the natural action of

${\bf PGL}(E) \times {\bf PGL}(F)$ on ${\bf P}^{8}$, which
does not change the rank of the $3 \times 3$ matrix $[z]$,
splits ${\bf P}^{8}$ into 3 orbits -- the fourfold $W$,
and the quasi-projective varieties $Sec(W) - W$
and ${\bf P}^{8} - Sec(W)$. The linear map ${\psi}$ sends
the closure of any of these orbits onto a determinantal
subvariety of ${\bf I}_{2}(W)$, and the only which has
to be seen is that these three determinantal subvarieties
are $D_{4}$,$D_{6}$ and $D_{9}$.

\smallskip
{\bf (6.3.3).}
{\sl
Quadrics of rank 6 related to the incidence correspondence
${\Sigma} \subset {\cal D} \times {\cal D}$.
}

Let $T$ be a general bidegree (2,2) divisor, and
let $Q \subset {\bf P}^{8}$ be any quadric such that
$Q \cap W = T$.  It follows from (6.3.2) that
such a quadric $Q$ is not unique --  $Q$ can be replaced
by any quadric in the projective 9-space
$Span(Q,{\bf I}_{2}(W))$,  i.e. $Q$ is unique
 mod.${\bf I}_{2}(W)$.

Let ${\Sigma} \subset G(7,9)$ (see (6.3.1)) be the incidence
correspondence

${\Sigma} = \ \mbox{(the closure of)}
\{ {\bf P}^{6} : C({\bf P}^{6}) = C + \tilde{C}, \  where \
 C, \tilde{C} \in {\cal D} \}$.
The set ${\Sigma}$ can be regarded (up-to closed subsets
 of codim.$ > 1$) also as an incidence correspondence
${\Sigma} \subset {\cal D} \times {\cal D}$.

Let ${\bf P}^{6} \in {\Sigma}$ be general.
The codim.2 subspace ${\bf P}^{6} \subset {\bf P}^{8}$
intersects the fourfold $W$ in an anticanonically
embedded del Pezzo surface $S({\bf P}^{6})$ of degree 6,
 and the quadric $Q$ intersects $S({\bf P}^{6})$
in a pair of elliptic sextics $C + \tilde{C}$.
Let ${\bf P}^{5}(C)$ and ${\bf P}^{5}(\tilde{C})$ be,
as in (6.3.1), the spans of $C$ and $\tilde{C}$, and
let $H$ and $\tilde{H}$ be linear forms on ${\bf P}^{8}$
such that $(H)_{o} \cap {\bf P}^{6} = {\bf P}^{5}(C)$,
$(\tilde{H})_{o} \cap {\bf P}^{6} = {\bf P}^{5}(\tilde{C})$.
The splitting $Q \cap S({\bf P}^{6}) = C + \tilde{C}$ implies:

$Q\mid_{{\bf P}^{6}} = H.{\tilde{H}}\mid_{{\bf P}^{6}}$
(mod. ${\bf I}_{2}(S({\bf P}^{6}))$ = the ``restriction'' of
${\bf I}_{2}(W)$ on ${\bf P}^{6}$).

(Here $Q$ is the quadratic form of $Q$, and we disregard
the multiplication by a non-zero constant).
Let ${\bf P}^{6} = (H_{1} = H_{2} = 0)$ be any pair of
linear equations which define the subspace
${\bf P}^{6} \subset {\bf P}^{8}$.
It follows from the preceding that $Q$ can be
represented in the form:

$Q = H.{\tilde{H}} + H_{1}.{\tilde{H}_{1}} + H_{2}.{\tilde{H}_{2}},
 mod.{\bf I}_{2}(W)$,
where $\tilde{H}_{1}$ and $\tilde{H}_{2}$ are some linear forms.
Clearly, the quadric
$H.{\tilde{H}} + H_{1}.{\tilde{H}_{1}} + H_{2}.{\tilde{H}_{2}}$
does not belong to the set $D_{6}(W)$ (the restriction
of this quadric to ${\bf P}^{6}$ does not contain the surface
$S({\bf P}^{6})$). In particular, the quadric  $Q$  in the definition
of $T$ can be replaced by this quadric of rank 6.

\smallskip
{\bf (6.4).}
{\sl
The Abel-Jacobi image $Z = {\Phi}({\cal D})$.
}
\smallskip

Let $T \subset W$, ${\cal D}$, etc., be as above, and
let ${\Phi}:{\cal D} \rightarrow J = J(T)$ be the Abel-Jacobi map
for the family ${\cal D}$.
Let $C \in {\cal D}$ be general, and let

${\Phi}^{*}:H^{1}(T,{\Omega}^{2}) \rightarrow
 H^{1}(N_{C/T} \otimes {\omega}_{T})$
be the codifferential of ${\Phi}$ in the
point $C \in {\cal D}$ (see e.g. [C]).

The space $H^{1}(T,{\Omega}^{2})$ is naturally isomorphic
to $H^{0}({\Omega}_{J(T)})$ -- the cotangent space
of $J(T)$ in a fixed point. The normal bundle
sequence for the embedding $T \subset W$, and the
formulae of Both and K\"unneth imply the isomorphism:

${\alpha}: H^{0}({\cal O}_{{\bf P}^{8}}(1)) \cong
 H^{0}(T,{\cal O}(1,1)) \rightarrow H^{1}(T,{\Omega}^{2})$.
(see also [Ve]). In particular, the elements of
$H^{1}(T,{\Omega}^{2})$ can be regarded as linear forms
on ${\bf P}^{8}$.

The following proposition is an analogue of Lemma 4.6 in [Vo]:

\smallskip
{\bf (6.4.1).}
{\bf Proposition.}
{\sl
Let $C \in {\cal D}$ be general, let ${\bf P}^{5}(C) = Span(C)$
be as in (6.3.1), let ${\bf P}^{6}$ be any 6-space through
${\bf P}^{5}(C)$, and let
$Q = H.{\tilde{H}} + H_{1}.{\tilde{H}_{1}} + H_{2}.{\tilde{H}_{2}},
mod.{\bf I}_{2}(W)$
be any of the representations of the quadric $Q$ defined by the
element ${\bf P}^{6} \in {\Sigma}$ (see (6.3.3)).

Let ${\Phi}^{*}$ and ${\alpha}$ be as above, and let
${\Phi}^{*}.{\alpha}$ be their composition.
Then the subspace
$ker({\Phi}^{*}.{\alpha}) \subset H^{0}({\cal O}_{{\bf P}^{8}}(1))$
is spanned on the forms
$H, \tilde{H}, H_{1}, \tilde{H}_{1}, H_{2}, \tilde{H}_{2}$.
}
\smallskip

{\bf Proof.}
see [Vo], the proof of Lemma 4.6. Note that the analogue of the
family ${\cal D}$, studied in [Vo], is the family of ``halves''
of canonical curves on the quartic double solid $B$  -- the family
of elliptic quartics on $B$ (see (6.3.1)).

\smallskip
{\bf (6.4.2).}
{\bf Corollary.}
{\sl
Let $Z = {\Phi}({\cal D})$ be the Abel-Jacobi image of the
family of elliptic sextics on $T$ of bidegree (3,3).
Then $dim(Z) = 3$.
}

\smallskip
{\bf (6.4.3).}
{\it
The fiber ${\Phi}^{-1}(z)$.
}

Let $C \in {\cal D}$ be general, let $z = {\Phi}(C)$, and let
${\Phi}^{-1}(z)_{o}$ be the connected component of
${\Phi}^{-1}(z)$ such that $C \in {\Phi}^{-1}(z)_{o}$.
Let ${\bf P}^{5}(C) = Span(C)$, and let $u,v$ be the local
parameters in the (general) point
${\bf P}^{7} = {\bf P}^{7}(0,0)
\in {\bf P}^{2}(C)^{*} = {\bf P}^{8}/{\bf P}^{5}(C)$.
Let $S(u,v) = T \cap {\bf P}^{7}(u,v)$. The surface
$S(u,v)$ is a hyperplane section of $T$ which contains
the (general) $C \in {\cal D}$. Such a surface $S(u,v)$
cannot be reducible. Otherwise
$l_{i} = p_{i}(C) = p_{i}(S(u,v)), i=1,2$
will be lines (see (6.2.1), (6.2.3)).
Therefore, $S(u,v)$ is a surface of type $K3$, and
the elliptic curve $C = C(u,v;0) \subset S(u,v)$ moves in
a pencil $\{ C(u,v;t) \subset S(u,v) \}$.
Obviously, $C(u,v;t) \in {\cal D}$ and
${\Phi}(C(u,v;t)) = {\Phi}(C)$,
since the curve $C(u,v;t)$ is a rational deformation
of $C = C(u,v;0)$. Therefore  $C(u,v;t)$ is the general
point of the irreducible 3-fold ${\Phi}^{-1}(z)_{o}$,
where $z = {\Phi}(C)$.

Let ${\bf P}^{6} \supset {\bf P}^{5}(C)$, and let
$Q_{o} = H.\tilde{H} + H_{1}.\tilde{H}_{1} + H_{2}.\tilde{H}_{2}$
be defined as in (6.3.3),(6.4); let $Q_{o}$ be also the
quadric defined by the equation $Q_{o} = 0$.
Let ${\Lambda}$ be the generator of $Q_{o}$ defined by the
condition ${\bf P}^{5}(C) \in {\Lambda}$. Obviously, the
5-space ${\bf P}^{5}(u,v;t) = Span(C(u,v,;t)) \in {\Lambda}$.
The correspondence
$C(u,v;t) \leftarrow {\bf P}^{5}(u,v;t)$
can be completed to a map
${\lambda}^{-1}:{\bf P}^{3}
\cong  {\Lambda} \rightarrow {\Phi}^{-1}(z)_{o}$.
It is easy to see that ${\lambda}^{-1}$  (hence ${\lambda}$)
is an isomorphism.

{\bf (*).}
It follows from the preceding that the quadric $Q_{o}$ does not
depend on the element $C(u,v;t) \in {\Phi}^{-1}(z)_{o}$.
We write $Q_{o} = Q_{o}(C) = Q_{o}(C(u,v;t)) = Q_{o}(z)$.

\smallskip
{\bf (6.5).}
{\bf Proposition.}
{\sl
Let $T \subset W$ be a general bidegree (2,2) divisor,
let $p_{i}:T \rightarrow {\bf P}^{2}, i=1,2,$  be the
projections, and let $Sing^{st}_{i}({\Theta})$ be as in (1.2.2).
Then:

  {\bf (i).}
There exist canonically defined maps
${\cal L}_{i}: Z \rightarrow {\cal L}_{i}(Z) \subset
 P(\tilde{\Delta}_{i},{\Delta}_{i}) \cong J(T)$,
where ${\cal L}_{i}(Z)$ is a component of
$Sing^{st}_{i}({\Theta})$, i=1,2.

  {\bf (ii).}
Let $C \in {\cal D}$ be general, and let $z = {\Phi}(C)$.
Then the quadric $Q_{o}(C) = Q_{o}(z)$ (see (6.4.3)(*))
coincides with the projectivized tangent cone
$Cone_{z}$ of ${\Theta}$ in the point
$z \in Z \subset Sing({\Theta})$.
}
\smallskip

{\bf Proof.}
Let $Q \subset {\bf P}^{8}$ be any quadric such that
$T = W \cap Q$, and let

${\bf I}_{2}(T) = {\bf P}( H^{0}({\bf P}^{8}, {\cal O}(2 - T)) )$
be the space of quadrics through $T$ (see also (6.3.2)). Clearly,

${\bf I}_{2}(T) \cong Span \{ {\bf I}_{2}(W), Q \}
\cong {\bf P}^{9}$.

Let
$D_{k}(T)
:= \ \mbox{(the closure of)} \{ P \in {\bf I}_{2}(T): rank(P) = k
\}$
be the $k$-th determinantal locus in  ${\bf I}_{2}(T)$.

Let $k = 6$. Then (see (6.3.2)(*))
 $D_{6}(T) \supset D_{6}(W) \cong Sec(W)$, and
$Q_{o}(z)$ does not belong to
${\bf I}_{2}(W) \supset D_{6}(W)$ (see (6.3.3)).
Therefore, the rule $z  \mapsto  Q_{o}(z)$ defines
a map $Q_{o}: Z \longrightarrow D_{6}(T)$,
and the image $Q_{o}(Z)$ is not a subset of
$D_{6}(W)$.

\smallskip
{\bf (6.5.1)}
{\bf Lemma.}
{\sl
$Q_{o}(Z)$ is a component of $D_{6}(T)$.
Moreover, $D_{6}(T) = Q_{o}(Z) \cup D_{6}(W)$.
}
\smallskip

{\bf Proof} of (6.5.1).
Let $Q_{o} \in D_{6}(T) - D_{6}(W)$, and let ${\Lambda}$
be one of the two generators of the quadric $Q_{o}$.
Let ${\bf P}^{5} \in {\Lambda}$. Then the set
$C({\bf P}^{5}) = T \cap {\bf P}^{5} =
 (W \cap Q) \cap {\bf P}^{5} =
 (W \cap Q_{o}) \cap {\bf P}^{5} =
 W \cap {\bf P}^{5}$
is a curve on $T$.  Moreover,  as it follows from the
elementary projective
properties of the fourfold  $W$,  $C({\bf P}^{5})$
is an elliptic curve of bidegree (3,3), i.e. $C \in {\cal D}$.
Clearly, $Q_{o} = Q_{o}(C({\bf P}^{5})) = Q_{o}(z)$,
where $z = {\Phi}(C({\bf P}^{5})) \in Z$. {\bf q.e.d.}

\smallskip
{\bf (6.5.2).}
{\it
The differential $dQ_{o}$ via the Gauss map of $Z$.
}

It follows from the definition of $Q_{o}(Z)$ that
$dim(Q_{o}(Z)) \le dim(Z) = 3$. Moreover, $Q_{o}(Z)$
is a component of the determinantal locus
$D_{6}(T) \subset {\bf I}_{2}(T) \cong {\bf P}^{9}$,
and the general quadric $Q \in {\bf I}_{2}(T)$ has
rank $9$. It follows from the general properties of
the determinantal varieties (see [F,ch.14]) that
the components of $D_{6}(T) = D_{9-3}(T)$ cannot be of
codimension greater than $3.(3+1)/2 = 6$. Therefore,
$dim(Q_{o}(Z)) = dim(Z) = 3$, and the map
$Q_{o}:Z \rightarrow Q_{o}(Z)$ is finite.
In particular, the differential
$dQ_{o}:T_{Z} \longrightarrow T_{Q_{o}(Z)}$
is a local isomorphism. Let
${\bf P}(dQ_{o}):{\bf P}(T_{Z})
\longrightarrow {\bf P}(T_{Q_{o}(Z)})$
be the projectivization of $dQ_{o}$.
Let

$Gauss:Z \longrightarrow G(3,9) = G(2:{\bf P}^{8})$,

$z \mapsto$  [ (the translate in $o \in J(T)$, of) the
set of all 3-spaces in the tangent space
of $Z \subset J(T)$ in the point $z \in Z $]

be the (rational) Gauss map, defined by the embedding
$Z \subset J(T)$.
Let $z \in Z$ be general. (In particular, $z$
is a regular point of $Z.$)
It follows from Proposition (6.4.1)
that the projective  2-space
${\bf P}^{2}(z) = vertex(Q_{o}(z))$
can be identified with
${\bf P}(Gauss(z)) \subset {\bf P}^{8} = {\bf P}(T_{J(T)}\mid_{o})$.
Since ${\bf P}(dQ_{o})$ is a local isomorphism,
we can identify the spaces ${\bf P}(T_{Q_{o}(Z)}\mid_{z})$ and
${\bf P}^{2}(z) = vertex(Q_{z})$ (see also [Ve, 4.29]).

\smallskip
{\bf (6.5.3).}
{\it
Proof of } (6.5)(i).

Let $i \in \{ 1,2 \}$ be fixed, and let
$L_{i}:{\cal D} \longrightarrow Symm^{18}(\tilde{\Delta}_{i})$
be the map defined in Lemma (6.2.3). (In (6.2.3), $p = p_{1}$ and
$L(C) = L_{1}(C)$;  the definition of $L_{2}(C)$ is evident.)
Let $C_{o} \in {\cal D}$ be general, and let $z = {\Phi}(C_{o})$.
Let ${\Lambda}$ be the generator of the quadric
$Q_{o}(z)$ defined by the condition
${\bf P}^{5}(C_{o}) \in {\Lambda}$.
Let

${\Lambda}_{i} = L_{i}({\Lambda})$ =
$\{ L_{i}(C): C = T \cap {\bf P}^{5}, {\bf P}^{5} \in {\Lambda} \}$
 $\subset Symm^{18}(\tilde{\Delta}_{i})$.

Clearly, the map $L_{i}:{\Lambda} \rightarrow {\Lambda}_{i}$
is an isomorphism. In particular,
${\Lambda}_{i} \cong {\bf P}^{3}$, i.e., the set
${\Lambda}_{i}$ is a rational subfamily of
$Supp \ {\Theta}_{i} \cup Supp \ P^{-}_{i}$ (see (6.2.3)(iii)). Here
${\Theta}_{i}$ and $P^{-}_{i}$ are the components of the
effective part in $Nm^{-1}({\omega}_{{\Delta}_{i}})$ -- see e.g. [W].
Therefore, all the divisors
$L_{i}(C) \in {\Lambda}_{i},
C = C({\bf P}^{5}), {\bf P}^{5} \in {\Lambda}$
belong to the same linear system
$\mid L_{i}(C_{o}) \mid$ on  $\tilde{\Delta}_{i}$.

Let
${\cal L}_{i} =
{\cal L}(L_{i}(C_{o})) \in {\bf Pic}^{18}(\tilde{\Delta}_{i})$
be the invertible sheaf defined by the effective divisor
$L_{i}(C_{o}) \in Symm^{18}(\tilde{\Delta}_{i})$.

The rule $L_{i}(C_{o}) \mapsto {\cal L}(L_{i}(C_{o}))$
defines a map

${\cal L}:L_{i}({\cal D}) \longrightarrow Supp \ {\Theta}_{i} \cup
Supp(P^{-})_{i}$.

It follows from Corollary (6.4.2), and from the
definition of the map ${\cal L} \circ L_{i}$ that
$dim({\cal L} \circ L_{i}({\cal D}) = dim(Z) = 3$. So, we obtained a
3-dimension
al
subset
${\cal L} \circ L_{i}({\cal D}) \ Nm^{-1}({\omega}_{{\Delta}_{i}})$
such that $h^{0}({\cal L}) \ge 4 $,  for any
${\cal L} \in {\cal L} \circ L_{i}({\cal D})$;
here we use the same symbol
${\cal L}$  for the sheaf  ${\cal L}$  and for the map  ${\cal L}.$

The number
$d = min \{ dim \mid {\cal L} \mid :
{\cal L} \in {\cal L} \circ L_{i}({\cal
D}) \}$  is a constant
throughout an open subset of ${\cal L} \circ L_{i}({\cal D})$.

\smallskip
{\bf (*).}
{\bf Lemma.}
{\sl
Let $d$ be as above. Then $d = 3$.
}
\smallskip

{\bf Proof} of (*).

Let, e.g., $d = 4$. (The case $d \ge 5$ can be treated
in a similar way.)
Let  $^{-}:\tilde{\Delta}_{i} \rightarrow \tilde{\Delta}_{i}$
be the involution induced by the double covering
$\tilde{\Delta}_{i} \rightarrow {\Delta}_{i}$,
and let

$W_{i} =
\{ {\cal M} \ = \ {\cal F} \otimes {\cal O}(x - \overline{x}):
 {\cal F} \in {\cal L} \circ L_{i}({\cal D}),
\  x \in \tilde{\Delta} \}$.

It follows from the definition of $W_{i}$ that
$W_{i} \subset Nm^{-1}({\omega}_{{\Delta}_{i}})$,
and $h^{0}({\cal M}) = 4$  for the general ${\cal M} \in W_{i}$
(see e.g. [Sh, Lemma 3.14]).

Therefore $W_{i}$ is a 4-dimensional subset of
$Sing^{st}_{i}({\Theta})$ (see (1.2.2)). However

$dim(Sing({\Theta})) = 3$ (see [Ve, Prop. 4.24]). Therefore
$d$ cannot be $4$. {\bf q.e.d.}

\smallskip
{\bf (**).}
{\bf Remark.}
The intermediate jacobian  $J(T)$ is a Prym variety which arises
from a double covering of a general plane sextic, in contrast to
the intermediate jacobian $J(B)$ of the desingularized nodal
quartic
double solid $B$ -- in which case the plane sextic $\Delta$
has a totally
tangent conic. The existence of a totally tangent conic
is a closed condition of codimension one, on the 19-dimensional
moduli space of the plane sextics.
Moreover,
 $dim(Sing({\Theta}(B))) = 4$ (see [De, 7]),
in contrast to $dim(Sing({\Theta}(T))) = 3.$
This, probably, once more explains why the Dixon correspondence,
which can be identified with a bidegree (2,2) divisor,
cannot be applied for a discriminantal pair which comes from
a nodal quartic double solid (see [Ve, 5]).

It follows from the preceding that the sheaf
${\cal L} \circ L_{i}(C)$
does not depend on the particular choice of the curve
$C \in {\Phi}^{-1}(z)_{o}$, $z = {\Phi}(C) = {\Phi}(C_{o})$.
Therefore, the map
${\cal L} \circ L_{i}:{\cal D} \longrightarrow J(T)$
factors
through the Abel-Jacobi map ${\Phi}:{\cal D} \rightarrow Z$.
Denote by ${\cal L}_{i}:Z \longrightarrow {\cal L}_{i}(Z) =
 {\cal L} \circ L_{i}({\cal D})$ the quotient map.
It follows from (*) that ${\cal L}_{i}(Z)$ is a 3-dimensional
component of $Sing^{st}_{i}({\Theta})$. This proves (i).

\smallskip
{\bf (6.5.4).}
{\bf
Proof } of (6.5)(ii).

Let ${\cal L}_{i}(z)$,
$z = {\Phi}(C_{o})$, etc., be as above. The sheaf
${\cal L}_{i}(z)$ is a stable singularity of $\Theta$,
with respect to $p_{i}$. Therefore, the projectivized
tangent cone
$Cone_{{\cal L}_{i}(z)}$ of $\Theta$, in the point
${\cal L}_{i}(z)$,
is a quadric which passes through the Prym-canonical image
${\Delta}_{i}^{T}$ of the discriminant curve ${\Delta}_{i}, i=1,2$
(see [Tju]).
Here we use the following results, due to Verra
(see [Ve, 3, and the proof
of 4.21]):

\smallskip
{\bf (*).}
{\sl
Let $s:{\Delta}_{i} \rightarrow T$ be the Steiner map, defined
by the rule:

$s: x  \mapsto  Sing(p_{i}^{-1}(x))$.

Then the image $s({\Delta}_{i})$ coincides with the Prym-canonical
curve ${\Delta}_{i}^{T}$.
}

{\bf (**).}
{\sl
Let $Q \subset {\bf P}^{8}$ be a quadric which passes through the
Steiner curves $s({\Delta}_{1})$ and $s({\Delta}_{2})$.
Then $Q \supset T$.
}
\smallskip

It follows from the preceding, and from (*) and (**), that

  {\bf (i).}
$Cone_{z} := Cone_{{\cal L}_{1}(z)} = Cone_{{\cal L}_{2}(z)}$;

  {\bf (ii).}
$Cone_{z} \supset T$,  i.e., $Cone_{z} \in {\bf I}_{2}(T)$.

It is well-known  that  $Cone_{z}$ is a quadric of rank 5 or 6
(see [K]), i.e. $Cone_{z} \in D_{5}(T) \cup D_{6}(T)$.
It is not hard to see that if $T$ is general then
$D_{5}(T) = \oslash$.  In fact, the general choice
of the quadric $Q$, such that $W \cap Q = T$, implies that
$codim(D_{5}(T) \subset {\bf I}_{2}(T)) = 10$, outside the
fixed determinantal $D_{4}(T) = D_{4}(W) \cong W$
(see (6.3.2), and [F, ch. 14]). Therefore,
$rank(Cone_{z}) = 6$.
The maps ${\cal L}_{1}$ and ${\cal L}_{2}$ are
local isomorphisms. Therefore the projective tangent spaces
$vertex(Q_{o}(z)) = {\bf P}^{2}(z) = {\bf P}(T_{Z}\mid_{z})$
(see (6.4.1)),  and
${\bf P}(T_{{\cal L}_{i}(Z)}\mid_{{\cal L}_{i}(z)}) =
 vertex(Cone_{{\cal L}_{i}(z)}) = vertex(Cone_{z}), i=1,2$
(see [Sh, 2.7 and 3.20]) can be identified.

It follows that the two quadrics:
$Cone_{z}$ and $Q_{o}(z)$
have the same vertex ${\bf P}^{2}(z)$, and
$Cone_{z} \supset T$, $Q_{o}(z) \supset T$. Moreover,
$Cone_{z}$ and $Q_{o}(z)$ belong to the determinantal
locus $D_{6}(T)$. An elementary projective consideration implies
that these two quadrics must coincide. This proves (ii).

As a corollary we obtain:

\newpage
{\bf (6.6)}
{\bf Theorem}
(The Torelli theorem for the Verra threefold).

{\sl
Let $T = T(2,2)$ be a general smooth bidegree (2,2) divisor in
the Segre image $W$ of ${\bf P}^{2} \times {\bf P}^{2}$, and
let $(J(T), {\Theta})$ be the principally polarized intermediate
jacobian of $T$. Then there exists a component
$Z$ of $Sing \ {\Theta}$ such that $dim(Z) = 3$, and
$T$ coincides with the intersection of all the projectivized
tangent cones of ${\Theta}$ in the regular points of $Z$.
}
\smallskip

{\bf Proof.}
It rests to be mentioned that
the map $Q_{o}$ sends the
set $Z = {\Phi}({\cal D})$ onto the component
$Q_{o}(Z)$ of $D_{6}(T)$
(see (6.3.3) and (6.5.1)),
and that the space ${\bf I}_{2}(T) \cong {\bf P}^{9}$
is spanned on the quadrics of the
determinantal locus $D_{6}(T).$
Moreover, $D_{6}(T) = Q_{o}(Z) \cup D_{6}(W)$,
and $Span(D_{6}(W))$ coincides with the proper subspace
${\bf I}_{2}(W) \cong {\bf P}^{8}$ of ${\bf I}_{2}(T)$
(see (6.3.2)(*) and (6.5.1)).  Then,  as it is not hard to see,
 $Span(Q_{o}(Z)) = {\bf I}_{2}(T)$.  Since the graded
ideal $I(T) = {\bigoplus} I_{d}(T)$ of $T \subset {\bf P}^{8}$
is generated by the component $I_{2}(T)$, and
${\bf I}_{2}(T) = {\bf P}(I_{2}(T))$,  the quadrics
of $Q_{o}(Z)$ (resp. -- the projective tangent cones
of ${\Theta}$ in the points of $Z$) cut the
projective subvariety $T \subset {\bf P}^{8}$ out
(see also [Vo, Prop.4.14]).

\smallskip
{\bf (6.7).}
{\bf Remarks.}

  {\bf (i).} (see [Ve, 4.17]):
Let
$Z_{T} = Sing^{st}_{1}({\Theta}) \cup Sing^{st}_{2}({\Theta})$.
Then $Z \subset Z_{T}$ can be separated among the components of
$Z_{T}$ by the numerical property:

{\sl
$Z$ = the union of all the irreducible components of $Z_{T}$
having not class $12.{\Theta}^{6}/6!$.
}

  {\bf (ii).}
Let $z \in Z$ be general, and let $C \in {\cal D}$ be such that
${\Phi}(C) = z$. Let
$Q_{o}(C) = Q_{o}(z) \in Q_{o}(Z) \subset D_{6}(T)$
be the rank 6 quadric attached to $z$, let ${\Lambda}$
be the generator of $Q_{o}(z)$ defined by
${\bf P}^{5}(C) \in {\Lambda}$, and let $\overline{\Lambda}$
be the complimentary generator of $Q_{o}(z)$. Let
$\overline{{\bf P}^{5}} \in \overline{\Lambda}$, let
$\overline{C} = T \cap \overline{{\bf P}^{5}}$, and let
$\overline{z} = {\Phi}(\overline{C})$ be the Abel-Jacobi image
of the curve $\overline{C} \in {\cal D}$. Obviously,
$Q_{o}(\overline{C}) = Q_{o}(\overline{z}) = Q_{o}(z)$, i.e.,
the degree of the finite map $Q_{o}:Z \rightarrow Q_{o}(Z)$ is
at least two. In fact, as it follows from the definition of the
map $Q_{o}$, the only preimages of the quadric $Q_{o}(z)$ are
the two points $z$ and $\overline{z}$, identified with the two
generators ${\Lambda}$ and $\overline{\Lambda}$ of the quadric.

\bigskip
\centerline{\sc 7.
The nodal $T(2,2)$.
}
\smallskip

{\bf (7.1).}
{\it
The tetragonal triples of Donagi connected with the
nodal $T(2,2)$.
}
\smallskip

Here we describe the two tetragonal triples which
correspond to the 4-gonal systems on the two nodal
discriminant sextics of the nodal $T(2,2)$
(see e.g. [Do]).

Let $T = W \cap Q$ has a simple node in the point
$(z)_{o} = (x)_{o} \times (y)_{o}$. Let
$p = p_{1}:T \rightarrow {\bf P}^{2}$
and
$q = p_{2}:T \rightarrow {\bf P}^{2}$
be the natural projections. Then the discriminant sextic
${\Delta}_{p}$ of $p$
(resp. -- ${\Delta}_{q}$  of $q$)
has a simple node in the point $(x)_{o}$
(resp. -- in the point $(y)_{o}$).
Let
${\bf P}^{1}_{p} = \mid{\cal O}_{{\bf P}^{2}}(1 - (x)_{o})\mid$
be the plane pencil of lines through $(x)_{o}$
(resp. -- ${\bf P}^{1}_{q} =
 \mid{\cal O}_{{\bf P}^{2}}(1 - (y)_{o})\mid$).
Let
$p^{-1}((x)_{o}) = L + \overline{L}$,
$q^{-1}((y)_{o}) = M + \overline{M}$.
Clearly,
$L \cap \overline{L} = M \cap \overline{M} = (z)_{o}$.

Let $pr:T \rightarrow {\bf P}^{7}$ be the rational
projection through $(z)_{o}$, and let $T_{+}$ be the
image of $T$. In particular, the proper images
$pr(L), pr(\overline{L}), pr(M)$, and $pr(\overline{M})$
are 4 isolated singular points of $T_{+}$ which lie on the
exceptional quadric $Q_{o} \subset T_{+}$ of  $pr$.

Let
$l \in {\bf P}^{1}_{p} - \{p(M),p(\overline{M})\}$,
$m \in {\bf P}^{2}_{q} - \{q(L),q(\overline{L})\}$,
and let $l$ and $m$ be, otherwise, {\it general}.
Let
$C(l,m) = p^{-1}(l) \cap q^{-1}(m) \cap T$.
A straightforward check gives that
 $C(l,m)$ is a curve of bidegree (2,2) and
of arithmetical genus one, which has a simple
node in the point $(z)_{o}$.
Let $q(l,m) = pr(C(l,m)) \subset T_{+}$
be the proper image of $C(l,m)$. It follows that
$q(l,m)$ is a conic. Thus, $T_{+}$ is birational
to a conic bundle
$s^{+}:T^{+} \rightarrow {\bf P}^{1}_{p} \times {\bf P}^{1}_{q}$.
It is not hard to describe the birational morphism
$T_{+} \rightarrow T^{+}$; it is a composition of
the blow-ups of the singular points $pr(L)$,
$pr(\overline{L})$,$pr(M)$,$pr(\overline{M})$, followed by
contracting of the four exceptional divisors along
their rulings. Because of the complexity of the notation,
caused by the additional exceptional sets, we shall
work on $T_{+}$, disregarding
the difference between $T_{+}$ and $T^{+}$; the statement will
not change substantially, if we work on $T^{+}$.

The (birational) conic bundle structure
$\{q(l,m): l,m \in {\bf P}^{1}_{p} \times {\bf P}^{1}_{q}\}$
on $T_{+}$ determines (the non-trivial component of)
the discriminant curve
${\Delta} \subset {\bf P}^{1}_{p} \times {\bf P}^{1}_{q}$
of  $s_{+}$.  Clearly, ${\Delta}$ is a smooth curve of
bidegree (4,4) on the quadric
${\bf P}^{1}_{p} \times {\bf P}^{1}_{q}$.

Let $l \in {\bf P}^{1}_{p}$ be general. Then
the surface
$S_{4}(l) = pr(p^{-1}(l)) \subset T_{+}$
is an anticanonically embedded del Pezzo
surface of degree 4. The map
$s_{+}:S_{4}(l) \rightarrow [l] \times {\bf P}^{1}_{q}$
defines a separated conic bundle structure on
$S_{4}(l)$, degenerated in the 4-tuple
${\Delta} \cap ([l] \times {\bf P}^{1}_{q})$.
Let $\tilde{\Delta}$ be the double discriminant curve
for $s_{+}$; $\tilde{\Delta}$ is isomorphic to the curve
of components of the degenerated fibers over $\Delta$.

Let
$g_{p} \subset Symm^{4}({\Delta})$
be the 4-gonal system on $\Delta$ defined by the
set of effective divisors
$\{ [l] \times {\bf P}^{1}_{q}: l \in {\bf P}^{1}_{p} \}$
(similarly -- for $g_{q}$),
and let
$s_{+}^{*}(g_{p}) =
 \{ L \in Symm^{4}(\tilde{\Delta}): (s_{+})_{*}(L) \in g_{p} \}.$
Let $l \in {\bf P}^{1}_{p}$, and $S_{4}(l)$ be as above.
The set of sixteen (-1)-curves on the anticanonically embedded
$S_{4}(l)$ coincide with the set of lines on $S_{4}(l)$.
The map $p = p_{1}$ defines a splitting of this set into two
``equal'' parts:  Eight of these curves come from the components
of the degenerated fibers
$p^{-1}(x), (x) \in (l \cap {\Delta}_{p}) - (x)_{o}$,
and the second 8-tuple is the set of these lines on  $S_{4}(l)$
which are
components of the four degenerated $s$-conics
$s_{+}^{-1}(u), u \in {\Delta} \cap ([l] \times {\bf P}^{1}_{q})$.
The lines from the second 8-tuple are components of fibers
of the conic bundle structure $s_{+}$. However, the lines
from the first 8-tuple are ``sections'' of $s_{+}$ -- the map
$s_{+}$ sends each of these lines isomorphically onto a line
on the base quadric ${\bf P}^{1}_{p} \times {\bf P}^{1}_{q}$.
Each of the lines of the first 8-tuple intersects
exactly 4 lines of the second 8-tuple. Moreover, if
${\lambda}_{1},{\lambda}_{2},{\lambda}_{3},{\lambda}_{4}$
is such a 4-tuple of lines (of the 2-nd system), then
$(s_{+})_{*}({\lambda}_{1} +...+ {\lambda}_{4}) =
 ([l] \times {\bf P}^{1}_{q}) \cap {\Delta}$.

The last causes a splitting of the natural preimage
$(s_{+})^{*}(([l] \times {\bf P}^{1}_{q}) \cap {\Delta})$
of $([l] \times {\bf P}^{1}_{q}) \cap {\Delta}$,
in $Symm^{4}(\tilde{\Delta})$,
into the following two subsets -- each of cardinality 8:

  {\bf (1).}  The set $\tilde{\Delta}^{+}_{p}(l)$ =
 the set of 4-tuples defined by the intersections
with the (projections of) the 8 components of the
degenerated fibers
$p^{-1}(x), (x) \in (l \cap {\Delta}_{p}) \ (x)_{o}$;

  {\bf (2).}  The complimentary set
$\tilde{\Delta}^{-}_{p}(l) :=
 (s_{+})^{*}(([l] \times {\bf P}^{1}_{q}) \cap {\Delta}) -
\tilde{\Delta}^{+}_{p}(l)$.

Clearly, this splitting does not depend on the
particular choice of the general line $l \in {\bf P}^{1}_{p}$.
Therefore, it defines a global splitting
$(s_{+})^{*}(g_{p}) =
 \tilde{\Delta}^{+}_{p} \cup \tilde{\Delta}^{-}_{p}$.

Obviously, the component
$\tilde{\Delta}^{+}_{p}$
is isomorphic to the non-singular model of the double
discriminant curve
$\tilde{\Delta}_{p}$
for the projection
$p:T \rightarrow {\bf P}^{2}$.

There are naturally defined involutions
$i^{+}_{p}:\tilde{\Delta}^{+}_{p} \rightarrow
 \tilde{\Delta}^{+}_{p}$ and
$i^{-}_{p}:\tilde{\Delta}^{-}_{p} \rightarrow
 \tilde{\Delta}^{-}_{p}$,
defined by interchanging of the 4-tuple
${\lambda}_{1} ,..., {\lambda}_{4}$
by its completion
$\overline{\lambda}_{1} ,..., \overline{\lambda}_{4}$.
(By definition,
${\lambda}_{i} + \overline{\lambda}_{i}, i = 1,...,4$,
are the four degenerated conics of
$s_{+}:S_{4}(l) \rightarrow ([l] \times {\bf P}^{1}_{q})$).

In fact, the 4-tuples
$({\lambda}_{1} ,..., {\lambda}_{4}) \in \tilde{\Delta}^{+}_{p}$
are (-1)-curves on $S_{4}(l)$;
the same -- for the complimentary 4-tuples.
Denote by
$S({\lambda}_{1},...,{\lambda}_{4})$ the ruled surface,
which is defined by contraction of the complimentary
4-tuple
$(\overline{\lambda}_{1} ,..., \overline{\lambda}_{4})$.
It follows from the definition of the elements of
$\tilde{\Delta}^{+}_{p}$ ( = the existence
of a secant (-1)-curve -- see above) that
$S({\lambda}_{1},...,{\lambda}_{4}) \cong {\bf F}_{1}$.

Similarly, the 4-tuples which belong to the component
$\tilde{\Delta}^{-}_{p}$
correspond to the relatively minimal models
$S( \overline{\lambda}_{1} ,..., \overline{\lambda}_{4} ),$
of the surfaces $S_{4}(l)$, which are of type
${\bf F}_{o}$ (i.e. -- quadrics).

Let
${\Delta}^{+}_{p} = \tilde{\Delta}^{+}_{p}/i^{+}_{p}$
and
${\Delta}^{-}_{p} = \tilde{\Delta}^{-}_{p}/i^{-}_{p}$
be the quotient curves.
Obviously, the natural 8-sheeted coverings
$\tilde{\Delta}^{+}_{p} \rightarrow {\bf P}^{1}_{p}$
and
$\tilde{\Delta}^{-}_{p} \rightarrow {\bf P}^{1}_{p}$
define the 4-sheeted coverings (the 4-gonal systems):
$g^{+}_{p}:{\Delta}^{+}_{p} \rightarrow {\bf P}^{1}_{p}$
and
$g^{-}_{p}:{\Delta}^{-}_{p} \rightarrow {\bf P}^{1}_{p}$.

In fact, we restored the tetragonal construction of
Donagi (see e.g.[Do]). Therefore, we proved the following
(see the notation above):

\smallskip
{\bf (7.2).}
{\bf Proposition.}

{\sl
$\{ (\tilde{\Delta},\Delta),
(\tilde{\Delta}^{+}_{p},{\Delta}^{+}_{p}),
(\tilde{\Delta}^{-}_{p},{\Delta}^{-}_{p}) \}$
is a 4-gonal triple of Donagi (see [Do]).
Moreover,
$\tilde{\Delta}^{+}_{p}$ is isomorphic to the
smooth model of the nodal discriminant plane sextic
${\Delta}_{p} \subset {\bf P}^{2}$, and the
involution
$i^{+}_{p}:\tilde{\Delta}^{+}_{p} \rightarrow
\tilde{\Delta}^{-}_{p}$
is a desingularization
of the involution $i_{p}$ on $\tilde{\Delta}_{p}$
(defined by the covering
$\tilde{\Delta}_{p} \rightarrow {\Delta}_{p}$).
}
\smallskip

Clearly, the same is true also for the 4-gonal
triple

$\{ (\tilde{\Delta}, \Delta),
(\tilde{\Delta}^{+}_{q},{\Delta}^{+}_{q}),
(\tilde{\Delta}^{-}_{q},{\Delta}^{-}_{q}) \}$
of Donagi, which corresponds to the
4-gonal system $g_{q}$ on the (4,4)-curve $\Delta$.
(Just like above, the curve $\tilde{\Delta}^{-}_{q}$
is isomorphic to the smooth model of the
double discriminant curve $\tilde{\Delta}_{q}$ of
$q:T \rightarrow {\bf P}^{2}$ ).

\smallskip
{\bf (7.3).}
{\bf Corollary.}
{\sl
Let $T = T(2,2)$ be a general nodal bidegree (2,2) divisor,
and let
$ (\tilde{\Delta}_{p},{\Delta}_{p})$
and
$ (\tilde{\Delta}_{q},{\Delta}_{q})$
be the discriminant pairs for the
natural projections
$p:T \rightarrow {\bf P}^{2}$
and
$q:T \rightarrow {\bf P}^{2}$.
Let
$ (\tilde{\Delta},{\Delta})$
be the discriminant pair of
the conic bundle structure
$s_{+}:T^{+} \rightarrow {\bf P}^{1} \times {\bf P}^{1}$
defined by the node of $T$ (see above),
and let $g_{p}$ and $g_{q}$ be the
4-gonal systems on the (4,4) curve $\Delta$
defined by the rulings of the quadric
${\bf P}^{1} \times {\bf P}^{1}$.
Then $g_{p}$ and $g_{q}$ define the two
tetragonal triples of Donagi:

$\{ (\tilde{\Delta},{\Delta}),
(\tilde{\Delta}^{+}_{p},{\Delta}^{+}_{p}),
(\tilde{\Delta}^{-}_{p},{\Delta}^{-}_{p}) \}$
and

$\{ (\tilde{\Delta},{\Delta}),
(\tilde{\Delta}^{+}_{q},{\Delta}^{+}_{q}),
(\tilde{\Delta}^{-}_{q},{\Delta}^{-}_{q}) \}$

such that
$ (\tilde{\Delta}^{+}_{p},{\Delta}^{+}_{p})$
is a desingularization of
the nodal pair
$ (\tilde{\Delta}_{p},{\Delta}_{p})$
(see above),
and the pair
$ (\tilde{\Delta}^{+}_{q},{\Delta}^{+}_{q})$
is a desingularization of
$ (\tilde{\Delta}_{q},{\Delta}_{q})$.
In other words, the Dixon correspondence
between the discriminant pairs of the
nodal bidegree (2,2) divisor $T$ is a composition
of two 4-gonal correspondences of Donagi.
}

\bigskip
{\bf Acknowledgment.}

I would like to express my gratitude to R.Donagi and
A.Verra for their substantial help.

\vspace{1.5in}

Institute of Mathematics

Bulgarian Academy of Sciences

ul. Acad. G.Bonchev, bl.8

1113 Sofia , Bulgaria

E-mail: ALGEBRA@BGEARN.BITNET

\end{document}